\documentclass[10pt]{article}
\usepackage{graphicx}
\usepackage{amsmath,amssymb,amsthm,amsfonts}
\usepackage{amssymb}
\newtheorem{thm}{Theorem}[section]

\newtheorem{lem}[thm]{Lemma}\newtheorem{prop}[thm]{Proposition}

\theoremstyle{definition}

\theoremstyle{remark}

\numberwithin{equation}{section}

\DeclareMathSymbol{\C}{\mathalpha}{AMSb}{"43}

\newcommand{\eps}{\varepsilon}

\newcommand{\lam}{\lambda}
\newcommand{\alp}{\alpha}

\newcommand{\R}{{\mathbb{R}}}
\newcommand{\N}{{\mathbb{N}}}

\def\ni{\noindent}

\def\qed{{\unskip\nobreak\hfil\penalty50
         \hskip2em\hbox{}\nobreak\hfil\mbox{\rule{1ex}{1ex} \qquad}
           \parfillskip=0pt
           \finalhyphendemerits=0\par }}
\def\proof{{\ni {\bf Proof:} } }

\newcommand{\bsub}{\begin{subequations}}
\newcommand{\esub}{\end{subequations}$\!$}
\begin{document}

\title{Compactness along the branch of semi-stable and unstable solutions for an elliptic problem with a singular nonlinearity \textsc \bf }

\author{\small \sc{Pierpaolo ESPOSITO\footnote{Dipartimento di Matematica,
Universit\`a degli Studi ``Roma Tre", Roma, Italy \& Pacific Institute for the Mathematical Sciences, University of British Columbia, Vancouver, B.C. Canada V6T 1Z2. E-mail: pesposito@pims.math.ca. Author supported by M.U.R.S.T., project
``Variational methods and nonlinear differential equations", and by a PIMS Postdoctoral Fellowship}\quad
Nassif GHOUSSOUB\footnote{Department of Mathematics, University of British Columbia, Vancouver, B.C. Canada V6T 1Z2. E-mail: nassif@math.ubc.ca. Research partially supported by the Natural Science and Engineering Research Council of Canada}\quad
Yujin GUO
\footnote{Department of Mathematics, University of British Columbia, Vancouver, B.C. Canada V6T 1Z2. E-mail: yjguo@math.ubc.ca. Author partially supported by the Natural Science Foundation of
P. R. China (10171036) and by a U.B.C. Graduate Fellowship.}}}
\date{\today}

\smallbreak
\maketitle

\begin{abstract}  \noindent We study the branch of semi-stable and unstable solutions (i.e., those whose Morse index is at most one)  of the Dirichlet boundary value problem  $-\Delta u=\frac{\lambda f(x)}{(1-u)^2}$ on a bounded domain $\Omega \subset \R^N$, which models --among other things-- a simple electrostatic Micro-Electromechanical System (MEMS) device. We extend the results of \cite{GG1} relating to the minimal branch, by obtaining compactness along unstable branches for $1\leq N \leq 7$ on any domain $\Omega$ and for a large class of ``permittivity profiles" $f$ .  We also show the remarkable fact that  power-like profiles $f(x) \simeq |x|^\alpha$ can push back  the critical dimension $N=7$ of this problem, by establishing compactness for the semi-stable branch on the unit ball, also for $N\geq 8$ and as long as $\alpha>\alpha_N=\frac{3N-14-4\sqrt{6}}{4+2\sqrt{6}}$ . As a byproduct, we are able to follow the second branch of the bifurcation diagram and prove the existence of a second solution for $\lambda$ in a natural range. In all these results, the conditions on the space-dimension and on the power of the profile are essentially sharp.
\end{abstract}

\vskip 0.2truein



\tableofcontents

\vskip 0.2truein

\section{Introduction}
We continue the analysis of \cite{GG1} for the problem:
$$
\left\{ \begin{array}{ll}
-\Delta u = \displaystyle\frac{\lambda
f(x)}{(1-u)^2}    & \hbox{in }\Omega ,\\
\quad \,\, 0<u<1 & \hbox{in } \Omega, \\
\quad \,\,\, u=0 & \hbox{on }\partial\Omega, \end{array} \right. \eqno{(S)_{\lam}}$$
where $\lambda>0$, $\Omega \subset \R^N$ is a bounded smooth domain and $f \in C(\bar \Omega)$ is a nonnegative function. This  equation models a
simple electrostatic Micro-Electromechanical System (MEMS) device
consisting of a thin dielectric elastic membrane with boundary
supported at $0$ below a rigid plate located at $+1$.   When
a voltage --represented here by $\lambda$-- is applied, the membrane
deflects towards the ceiling plate and a snap-through may occur when
it exceeds a certain critical value $\lambda^*$ (pull-in voltage).
This creates a so-called  ``pull-in instability" which greatly
affects the design of many devices (see \cite{FMP,PB} for a detailed discussion on MEMS devices). The mathematical model lends to
a nonlinear parabolic problem for the dynamic deflection of the
elastic membrane which has been considered by the second and third-named authors in \cite{GG2,GG3}. Concerning the stationary equation, in \cite{GG1} the challenge was to estimate $\lambda^*$ in terms of
material properties of the membrane, which can be fabricated with a
spatially varying dielectric permittivity profile $f(x)$. In particular, lower bounds for $\lam^*$ were proved completing in this way the upper bounds of \cite{GPW,P1}. In all the above-mentioned papers, one can recognize  a clear distinction --in techniques and in the available results-- between the case where the permittivity profile $f$ is bounded away from zero, from where it is allowed to vanish somewhere. A test case for the latter situation --that has generated much interest among both mathematicians and engineers-- is when we have a  power-law permittivity profile $f(x)=|x|^\alpha$ ($\alpha \geq 0$) on a ball.

\medskip \noindent There already exist in the litterature many interesting results concerning the properties of the branch of semi-stable solutions for Dirichlet boundary value problems of the form $-\Delta u = \lambda h(u)$ where $h$ is a regular nonlinearity (for example of the form $e^u$ or $(1+u)^p$ for $p>1$). See for example the seminal papers  \cite{CR,JL,KK} and also \cite{C} for a survey on the subject and an exhaustive list of related references. The singular situation was considered in a very general context in \cite{MP}, and this analysis was completed  in \cite{GG1}  to allow for a general continuous permittivity profile $f(x) \geq 0$. Fine properties of steady states --such
as regularity, stability, uniqueness, multiplicity, energy estimates
and comparison results-- were shown there to depend on the dimension of the
ambient space and on the permittivity profile.

\medskip \noindent Let us fix some notations and terminology. The  {\it minimal solutions} of the equation are those classical solutions $u_{\lam}$ of $(S)_{\lam}$ that satisfy $u_{\lam}(x)\le u(x)$ in $\Omega $ for any solution $u$ of $(S)_{\lam}$.  Throughout and unless otherwise specified,  solutions for $(S)_\lambda$ are considered to be in the classical sense. Now for any solution $u$ of $(S)_{\lam}$, one can introduce the linearized operator at $u$ defined by:
$$L_{u, \lambda} =-\Delta -\frac{2\lambda f(x)}{(1-u)^3},$$
and its corresponding eigenvalues $\{\mu_{k, \lambda}(u); k=1, 2,...\}$. Note that  the first eigenvalue is simple and is given by:
\[
\mu_{1,\lambda}(u)= \inf \left\{ \left\langle L_{u,\lambda} \phi,\phi \right\rangle_{H_0^1(\Omega)}; \, \phi \in C_0^\infty(\Omega),  \int_\Omega  |\phi (x)|^2dx =1\right\}
\]
with the infimum being attained at a first eigenfunction $\phi_1$, while  the second eigenvalue is given by the formula:
\[
\mu_{2,\lambda}(u)=
  \inf \left\{ \left\langle L_{u,\lambda} \phi,\phi \right\rangle_{H_0^1(\Omega)};  \, \phi \in C_0^\infty(\Omega),  \int_\Omega  |\phi (x)|^2 dx=1\, {\rm and}\, \int_\Omega \phi (x) \phi_1 (x)dx=0\right\}.
  \]
This construction can then be iterated to obtain the $k$-th eigenvalue $\mu_{k,\lambda}(u)$ with the convention that eigenvalues are repeated according to their multiplicities.

\medskip  \noindent The usual analysis of the minimal branch (composed of semi-stable solutions) was extended in
\cite{GG1} by Ghoussoub and Guo to cover the singular situation
$(S)_\lambda$ above and  the subsequent result -- best illustrated
by the following bifurcation diagram-- was obtained.

\begin{figure}[htbp]
\begin{center}
{\includegraphics[width = 8cm,height=5.5cm,clip]{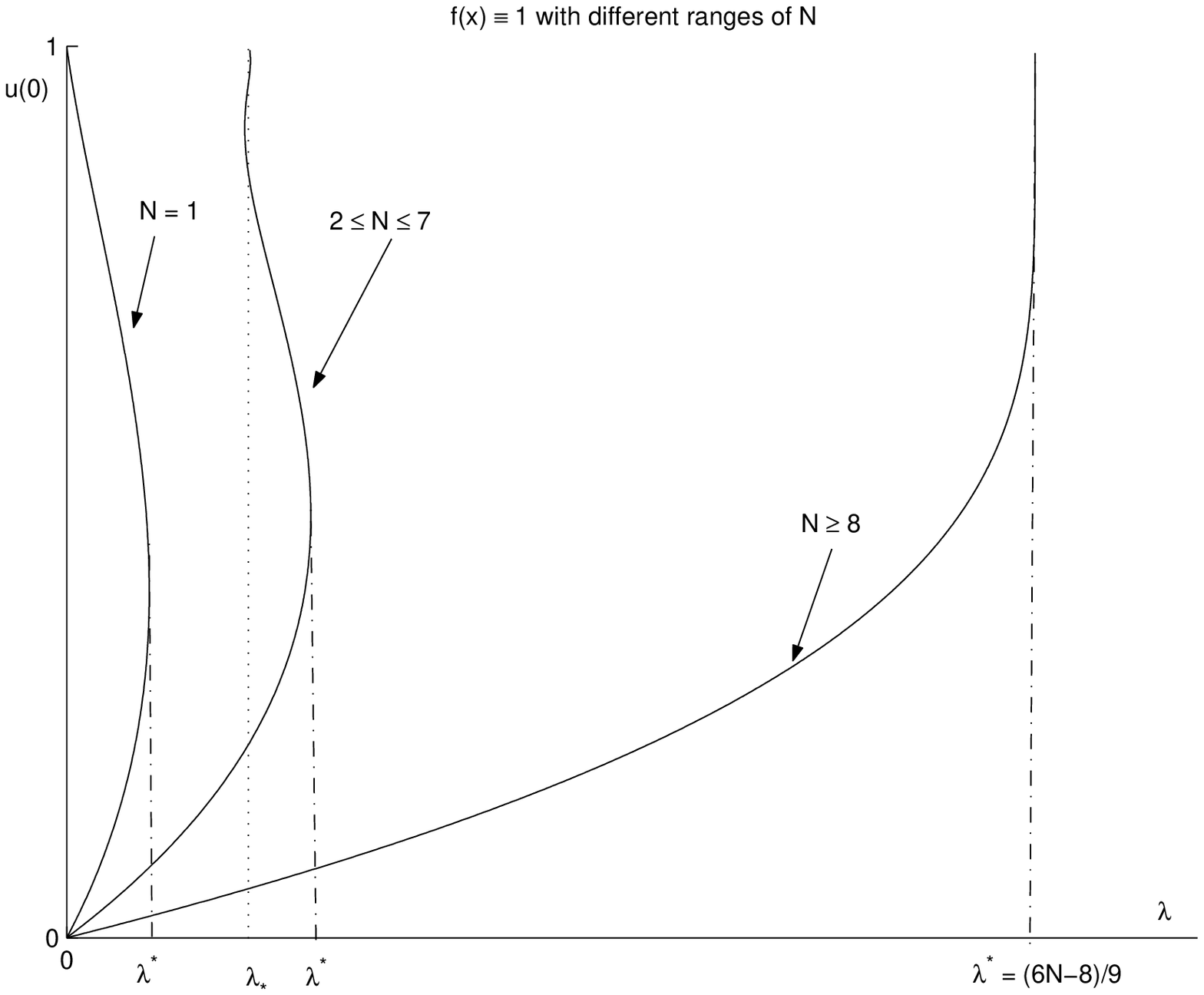}
\label{fig:fig2_2a}}
{\includegraphics[width = 8cm,height=5.5cm,clip]{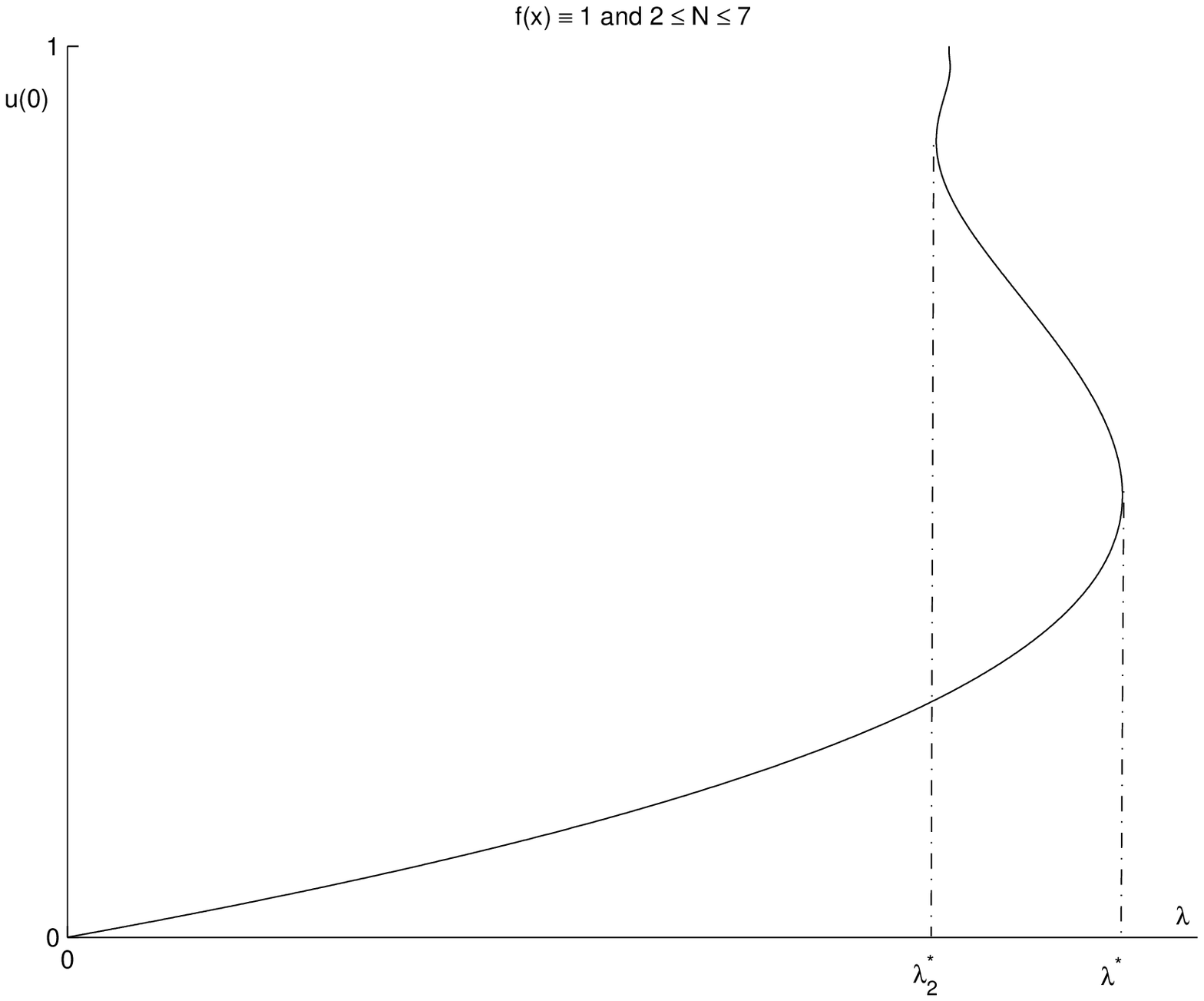}
\label{fig:fig2_2b}} \caption{{\em Left figure: plots of $u(0)$
versus $\lam$ for the case where $f(x)\equiv 1$ is defined in the
unit ball $B_1(0)\subset \R^N$ with different ranges of dimension
$N$, where we have $\lam ^*=(6N-8)/9$ for dimension $N\ge 8$. Right
figure: plots of $u(0)$ versus $\lam$ for the case where $f(x)\equiv
1$ is defined in the unit ball $B_1(0)\subset \R^N$ with dimension
$2\le N\le 7$, where $\lam ^*$ (resp.~$\lam ^*_2$) is the first
(resp.~second) turning point.}} \label{fig:fig2_2}
\end{center}
\end{figure}

\noindent {\bf Theorem A (Theorem 1.1-1.3 in \cite{GG1}): } {\em Suppose $f \in C(\bar \Omega)$ is a nonnegative
function on $\Omega$. Then, there exists a
finite $\lam ^*>0$ such that
\begin{enumerate}
\item If $0\leq \lam <\lam ^*$, there exists a unique minimal solution $u_\lambda$ of
$(S)_{\lam}$ such that $\mu_{1,\lambda}(u_\lambda)>0$;
\item If $\lam >\lam ^*$, there is no solution for $(S)_{\lam}$.

\item Moreover, if $1\leq N\leq 7$ then --by means of energy estimates-- one has
\begin{equation}
\label{estim} \sup_{\lambda \in (0,\lambda^*)} \parallel u_\lambda \parallel_\infty <1   \end{equation}
and consequently,   $u^*=\displaystyle\lim_{\lambda \uparrow \lambda^*}u_\lambda$ is a solution for $(S)_{\lam^*}$ such that
\begin{equation}
\mu_{1,\lambda^*}(u^*)=0.  \label{minimal*}
\end{equation}
In particular,   $u^*$  --often referred to as the extremal solution of problem $(S)_{\lam }$-- is unique.
\item On the other hand, if $f(x)=|x|^\alpha$ and $\Omega$ is the unit ball, then  the extremal solution is necessarily
$u^*(x)=1-|x|^{\frac{2+\alpha}{3}}$ and $\lambda^*=\frac{(2+\alpha)(3N+\alpha-4)}{9}$, provided $N\geq 8$ and $0\leq \alpha \leq \alpha_N=\frac{3N-14-4\sqrt{6}}{4+2\sqrt{6}}$.
\end{enumerate}}

\medskip \noindent We note that in general, the function $u^*$ exists in any dimension, does solve $(S)_{\lam^*}$ in an appropriate weak sense and is the unique solution in some suitable class (see the Appendix).

\medskip \noindent Our first goal is the study of the effect of power-like permittivity profiles $f(x) \simeq |x|^\alpha$ for the problem $(S)_\lam$ on the unit ball $B=B_1(0)$. We extend the previous result in higher dimensions:
\begin{thm} \label{thm0} Assume $N \geq 8$ and $\alpha>\alpha_N=\frac{3N-14-4\sqrt{6}}{4+2\sqrt{6}}$. Let $f \in C(\bar B)$ be such that:
\begin{equation} \label{assf0}  f(x)=|x|^{\alpha} g(x)\:,\quad g(x) \geq C>0 \hbox{ in }B. \end{equation}
Let $(\lambda_n)_n$ be such that $\lambda_n \to  \lam \in [0,\lam^*]$ and $u_n$ be a solution of $(S)_{\lam_n}$ so that $\mu_{1,n}:=\mu_{1,\lam_n}(u_n) \geq 0$. Then,
$$\sup_{n \in \N} \parallel u_n \parallel_\infty <1.$$
In particular, the extremal solution $u^*=\displaystyle\lim_{\lam \uparrow \lam^*} u_\lam$ is a solution of $(S)_{\lam^*}$ such that (\ref{minimal*}) holds.
\end{thm}

\medskip \noindent As to non-minimal solutions, it is also shown in \cite{GG1} --following ideas of Crandall-Rabinowitz \cite{CR}-- that, for  $1 \leq N \leq 7$, and  for $\lambda$ close enough to $\lambda^*$, there exists a unique second branch $U_\lambda$ of solutions for $(S)_{\lam}$, bifurcating from $u^*$, such that
\begin{equation}  \mu_{1,\lambda}(U_\lambda) < 0  \quad {\rm while} \quad \mu_{2,\lambda}(U_\lambda)> 0.  \label{nonminimalev}
\end{equation}
For $N\geq 8$ and $\alpha >\alpha_N$, the same remains true for problem $(S)_{\lam}$ on the unit ball with $f(x)$ as in (\ref{assf0}) and $U_\lam$ is a radial function.

\medskip \noindent In the sequel, we try to provide a rigorous analysis for other features of the bifurcation diagram,  in particular the second branch of unstable solutions, as well as the second bifurcation point. But first, and for  the sake of completeness, we shall give a variational characterization for the unstable solution $U_\lambda$ in the following sense:
\begin{thm} \label{thm1} Assume $f$ is a non-negative function in $C(\bar\Omega)$ where $\Omega$ is a bounded domain in $\R^N$. If $1\le N\le 7$, then there exists $\delta>0$ such that  for any $\lambda \in (\lambda^*-\delta,\lambda^*)$,  the second solution $U_\lambda$ is a mountain pass solution for some regularized energy functional $J_{\varepsilon,\lambda}$ on the space $H_0^1(\Omega)$.

Moreover,  the same result is still true for $N\geq 8$ provided $\Omega$ is a ball, and $f(x)$ is as in (\ref{assf0}) with $\alpha>\alpha_N$.
\end{thm}

\noindent We are now interested in continuing  the second branch till the second bifurcation point, by means of the Implicit Function Theorem. For that, we have the following compactness result:
\begin{thm} \label{thm2} Assume $2\leq N \leq 7$. Let $f \in C(\bar \Omega)$ be such that:
\begin{equation} \label{assf1}  f(x)=\left(\prod_{i=1}^k |x-p_i|^{\alpha_i}\right) g(x)\:,\quad g(x) \geq C>0 \hbox{ in }\Omega, \end{equation}
for some points $p_i \in \Omega$ and exponents $\alpha_i \geq 0$. Let $(\lambda_n)_n$ be a sequence such that $\lambda_n \to  \lam \in[0,\lam^*]$ and let $u_n$ be an associated solution such that
\begin{equation} \label{morse}
\mu_{2,n}:=\mu_{2,\lambda_n}(u_n)\geq 0.  \end{equation}
Then, $\displaystyle\sup_{n \in \N} \parallel u_n \parallel_\infty <1$. Moreover, if in addition $\mu_{1,n}:=\mu_{1,\lambda_n}(u_n)<0$, then necessarily $\lambda >0$.
\end{thm}
\noindent Let us mention that Theorem \ref{thm2} yields another proof --based on a blow-up argument-- of the compactness result for minimal solutions  (\ref{estim}) established in \cite{GG1} by means of some energy estimates, though under the more stringent assumption (\ref{assf1}) on $f(x)$. We expect that the same result should be true for radial solutions on the unit ball for $N\geq 8$, $\alpha>\alpha_N$, and $f \in C(\bar \Omega)$ as in (\ref{assf0}).

\medskip \noindent As far as we know, there are no compactness results of this type in the case of regular nonlinearities,  marking a substantial difference with  the singular situation. Theorem \ref{thm2} is based on a blow up argument and the knowledge of linear instability for solutions of a limit problem on $\R^N$, a result which is interesting in itself (see for example \cite{CC}) and which somehow explains the special role of dimension $7$ and $\alpha=\alpha_N$ for this problem.
\begin{thm} \label{thm3} Assume that either $1\le N\le 7$ and $\alpha \geq 0$ or that $N\geq 8$ and $\alpha>\alpha_N$. Let $U$ be a solution of
\begin{equation} \left\{ \begin{array}{ll} \Delta U = \displaystyle\frac{|y|^\alpha}{U^2}
& \hbox{in }\R ^N,\\
U(y) \geq C>0  & \hbox{in } \R^N. \end{array} \right.  \label{3:1agg}
\end{equation}
Then, $U$ is linearly unstable in the following sense:
\begin{equation}
\mu _1(U)=\inf \left\{ \int_{\R^N} \big( |\nabla \phi
|^2- \frac{2 |y|^\alpha}{U^3} \phi ^2\big) dx; \, \phi \in C^{\infty }_0(\R^N)\, {\rm and}\,  \int_{\R^N} \phi^2=1 \right\}<0\,. \label{3:2agg}
\end{equation}
Moreover, if $N\geq 8$ and $0\leq \alpha\leq \alpha_N$, then there exists at least a  solution $U$ of (\ref{3:1agg}) such that $\mu_1(U) \geq 0$.
\end{thm}

\medskip \noindent Theorem \ref{thm3} is the main tool to control the blow up behavior of a possible non compact sequence of solutions. The usual asymptotic analysis for equations with Sobolev critical nonlinearity, based on some energy bounds (usually $L^{\frac{2N}{N-2}}(\Omega)$-bounds), does not work in our context. In view of \cite{GG1}, a possible loss of compactness can be related to the $L^{\frac{3N}{2}}(\Omega)$-norm along the sequence. Essentially, the blow up associated to a sequence $u_n$ (in the sense of the blowing up of $(1-u_n)^{-1}$) corresponds exactly to the blow up of the $L^{\frac{3N}{2}}(\Omega)$-norm. We replace these energy bounds by some spectral information and, based on Theorem \ref{thm3}, we provide an estimate of the number of blow up points (counted with their ``multiplicities'') in terms of the Morse index along the sequence.

\medskip \noindent We now define the second bifurcation point in the  following way for $(S)_\lam$:
$$\lambda_2^*=\inf \{ \beta>0: \exists \hbox{ a curve} \, V_\lambda \in C([\beta ,\lambda^*]; C^2(\Omega)) \hbox{ of solutions for } (S)_\lambda \hbox{ s.t. }\mu_{2,\lambda}(V_\lambda) \geq 0  , \:V_\lambda \equiv U_\lambda \, \forall \lambda \in (\lambda^*-\delta,\lambda^*) \}.$$
We then have the following multiplicity result:
\begin{thm} \label{thm4}
Assume $f \in C(\bar \Omega)$ to be of the form (\ref{assf1}). Then, for $2\leq N \leq 7$ we have that  $\lambda_2^* \in (0,\lambda^*)$ and for any $\lambda \in (\lambda_2^*,\lambda^*)$ there exist at least two solutions $u_\lambda$ and  $V_\lambda$ for $(S)_\lambda$, so that
$$\mu_{1,\lam}(V_\lam)<0 \quad {\rm while}\quad  \mu_{2,\lam}(V_\lam)\geq 0.$$
In particular, for $\lambda=\lambda_2^*$, there exists a second solution, namely $V^*:=\displaystyle\lim_{\lambda \downarrow \lambda_2^*} V_\lambda$ so that
$$\mu_{1,\lambda_2^*}(V^*)<0  \quad {\rm and}\quad \mu_{2,\lambda_2^*}(V^*)=0.$$
\end{thm}

\medskip \noindent One can compare Theorem \ref{thm4} with the multiplicity result of \cite{ABC} for nonlinearities of the form $\lambda u^q+u^p$ ($0<q<1<p$), where the authors show that  for $p$ subcritical, there exists a second --mountain pass-- solution for any $\lam \in [0,\lam^*)$. On the other hand, when $p$ is critical, the second branch blows up as $\lambda \to 0$ (see also \cite{AT1} for a related problem). We note that in our situation, the second branch cannot approach the value $\lam=0$ as illustrated by the bifurcation diagram above.

\medskip \noindent Let now $V_\lambda$, $\lambda \in (\beta, \lambda^*)$ be one of the curves appearing in the definition of $\lambda_2^*$. By (\ref{nonminimalev}), we have that $L_{V_\lam,\lam}$ is invertible for $\lam \in (\lam^*-\delta,\lam^*)$ and, as long as it remains invertible, we can use the Implicit Function Theorem to find $V_\lam$ as the unique smooth extension of the curve $U_\lambda$ (in principle $U_\lambda$ exists only for $\lambda$ close to $\lambda^*$). Let now $\lam^{**}$ be defined in the following way:
$$\lambda^{**}=\inf \{ \beta>0:\forall \lambda \in (\beta ,\lambda^*) \: \exists \: V_\lambda \hbox{ solution of } (S)_\lambda \hbox{ so that }\mu_{2,\lambda}(V_\lambda) >0  , \:V_\lambda \equiv U_\lambda \, {\rm for}\, \lambda \in (\lambda^*-\delta,\lambda^*) \}.$$
Then, $\lam_2^* \leq \lam^{**}$ and there exists a smooth curve $V_\lam$ for $\lam \in (\lam^{**},\lam^*)$ so that $V_\lam$ is the unique maximal extension of the curve $U_\lambda$. This is what the second branch  is supposed to be. If now $\lam_2^*<\lam^{**}$, then for $\lam \in (\lam_2^*,\lam^{**})$ there is no longer uniqueness for the extension and the ``second branch'' is defined only as one of potentially many continuous extensions of $U_\lam$.

\medskip \noindent It remains open the problem  whether $\lambda_2^*$ is the second turning point for the solution diagram of $(S)_\lambda$ or if the ``second branch'' simply disappears at $\lambda=\lambda_2^*$. Note that if the ``second branch'' does not disappear, then it can continue for $\lam$ less than $\lam_2^*$ but only along solutions whose first two eigenvalues are negative.

\medskip \noindent In dimension $1$, we have a stronger but somewhat different compactness result. Recall that  $\mu_{k,\lambda_n}(u_n)$ is the $k-$th eigenvalue of $L_{u_n,\lambda_n}$ counted with their multiplicity.

\begin{thm} \label{thm5}Let $I$ be a bounded interval in $\R$ and $f \in C^1(\bar I)$ be such that $f \geq C>0$ in $I$.
Let $(u_n)_n$ be a solution sequence for $(S)_{\lambda_n}$ on $I$, where $\lambda_n \to  \lam \in[0,\lam ^*]$. Assume that for any $n \in \N$ and $k$ large enough, we have:
\begin{equation} \label{morsebound}
\mu_{k,n}:=\mu_{k,\lambda_n}(u_n) \geq 0.
\end{equation}
If $\lambda >0$, then again  $\displaystyle\sup_{n \in \N} \parallel u_n \parallel_\infty <1$ and compactness holds.
\end{thm}

\medskip \noindent Even in dimension $1$, we can still define $\lambda_2^*$ but we don't know when $\lambda_2^*=0$  (this is indeed the case when $f(x)=1$, see \cite{Z}) or when $\lambda_2^*>0$. In the latter situation, there would exist a solution $V^*$ for $(S)_{\lambda_2^*}$ which could be --in some cases-- the second turning point. Let us remark that the multiplicity result of Theorem \ref{thm4} holds also in dimension $1$ for any $\lambda \in (\lambda_2^*,\lambda^*)$.

\medskip \noindent The paper is organized as follows. In Section $2$ we
provide the mountain pass variational characterization of $U_\lambda$ for $\lambda$ close to $\lambda^*$ as stated in Theorem \ref{thm1}. The compactness result of Theorem \ref{thm0} on the unit ball is proved in Section $3$. Section $4$ is concerned with the compactness of the second branch of $(S)_\lam$ as stated in Theorem \ref{thm2}. Section $5$ deals with the dimension $1$ of Theorem \ref{thm5}. In Section $6$ we give the proof of the multiplicity result in Theorem \ref{thm4}. Finally, the linear instability property of Theorem \ref{thm3} and the details of the above mentioned counterexample to the $C^2$-regularity of $u^*$ in dimension $N\geq 8$, $0\leq \alpha \leq \alpha_N$, are given in the Appendix.

\section{Mountain Pass solutions}
This Section is devoted to the variational characterization of the second solution $U_\lambda$ of $(S)_{\lam }$ for $\lambda \uparrow \lambda^*$ and in dimension $1\le N\le 7$. Let us stress that the argument works also for problem $(S)_\lam$ on the unit ball with $f(x)$ in the form (\ref{assf0}) provided  $N\geq 8$, $\alpha>\alpha_N$.

\medskip \noindent Since the nonlinearity $g(u)=\frac{1}{(1-u)^2}$ is singular at $u=1$, we need to consider a regularized $C^1$ nonlinearity $g_\varepsilon(u)$, $0<\varepsilon<1$, of the following form:
\begin{equation} \label{geps} g_\varepsilon(u)=\left\{ \begin{array}{ll}
\frac{1}{(1-u)^2} & u\le 1-\varepsilon \,,\\
\frac{1}{\varepsilon ^2}-\frac{2(1-\varepsilon )}{p\varepsilon
^3}+\frac{2}{p\varepsilon ^3(1-\varepsilon )^{p-1}}u^p &
u\ge 1-\varepsilon\,,\end{array} \right. \end{equation}
where $p>1$ if $N=1,2$ and $1<p<\frac{N+2}{N-2}$ if $3\leq N \leq 7$. For $\lam\in (0, \lam ^*)$, we study the regularized semilinear elliptic problem:
\begin{equation}\left\{ \begin{array}{ll}-\Delta u =\lam f(x)  g_\varepsilon(u) & \hbox{in }\Omega , \\
\quad \,\,\, u=0  & \hbox{on } \partial\Omega. \end{array} \right. \label{2:1} \end{equation}
From a variational viewpoint, the action functional associated to (\ref{2:1}) is
\begin{equation}
J_{\varepsilon, \lam}(u)=\displaystyle
\frac{1}{2}\int _{\Omega}|\nabla u|^2dx -\lam \int _{\Omega}f(x) G_{\varepsilon }(u)dx\,,\quad
u\in H^1_0(\Omega )\,, \label{2:2}
\end{equation}
where $G_{\varepsilon }(u)=\displaystyle\int ^u_{-\infty}\displaystyle
g_{_\varepsilon }(s)ds$.

\medskip \noindent Fix now $0<\varepsilon <\frac{1-\parallel u^*\parallel_\infty}{2}$. For $\lam \uparrow \lam^*$, the minimal solution $u_\lam $ of $(S)_\lam$ is still a solution of (\ref{2:1}) so that $\mu_1\left( -\Delta-\lambda f(x) g_\varepsilon'(u_\lambda) \right)>0$. In order to motivate the choice of $g_\varepsilon(u)$, we briefly sketch the proof of Theorem \ref{thm1}. First, we prove that $u_\lambda$ is a local minimum for $J_{\varepsilon,\lambda}(u)$ for $\lambda \uparrow \lambda^*$. Then, by the well known Mountain Pass Theorem \cite{AR}, we show the existence of a second solution $U_{\varepsilon,\lambda}$ for (\ref{2:1}). Since $U_{\varepsilon,\lambda} \to u^*$ in $C(\bar \Omega)$ as $\lambda \uparrow \lambda^*$, we have that $U_{\varepsilon,\lambda}\leq 1-\varepsilon$ and  $U_{\varepsilon,\lambda}$ is then a second solution for $(S)_\lambda$ bifurcating from $u^*$. But since $U_{\varepsilon,\lambda}$ is a MP solution and $(S)_\lambda$ has exactly two solutions $u_\lambda$, $U_\lambda$ for $\lambda \uparrow \lam^*$, we get that $U_{\varepsilon,\lambda}=U_\lambda$.\\
The subcritical growth:
\begin{equation} 0 \leq g_\varepsilon(u)\leq C_\varepsilon(1+|u|^p)\label{growth}\end{equation}
and the inequality:
\begin{equation} \theta G_\varepsilon(u) \leq u g_\varepsilon(u) \quad \hbox{for }u \geq M_\varepsilon, \label{cruci}\end{equation}
for some $C_\varepsilon,\: M_\varepsilon>0$ large and $\theta=\frac{p+3}{2}>2$,  will yield that $J_{\varepsilon,\lambda}$ satisfies the Palais-Smale condition and, by means of a bootstrap argument, we get the uniform convergence of $U_{\varepsilon,\lambda}$. On the other hand, the convexity of $g_\varepsilon(u)$ ensures that problem (\ref{2:1}) has the unique solution $u^*$ at $\lambda=\lambda^*$, which then allows us to identify the limit of $U_{\varepsilon,\lambda}$ as $\lam \uparrow \lam^*$.

\medskip \noindent In order to complete the details for the proof of Theorem \ref{thm1}, we first need to show the following:
\begin{lem} For $\lam \uparrow \lam ^*$, the minimal solution $u_{\lam}$ of
$(S)_{\lam}$ is a local minimum of $J_{\varepsilon,\lambda}$ on $H^1_0(\Omega)$.
\end{lem}
\proof First, we show that $u_\lambda$ is a local minimum of $J_{\varepsilon,\lambda}$ in $C^1(\bar \Omega)$. Indeed,  since $$\mu_{1,\lambda}:=\mu_1 \left( -\Delta-\lambda f(x) g_\varepsilon'(u_\lambda) \right)>0,$$
we have the following inequality:
\begin{equation}
\int _{\Omega}|\nabla \phi|^2 dx-2\lam \int _{\Omega} \frac{f(x)}{(1-u_\lam)^3} \phi ^2 dx\ge \mu_{1,\lam}\int_\Omega \phi^2
\label{2:4} \end{equation}
for any $\phi \in H^1_0(\Omega)$, since $u_\lam \leq 1-\eps$. Now, take any $\phi \in H_0^1(\Omega) \cap C^1(\bar \Omega)$ such that $\| \phi \|_{C^1} \leq \delta_\lam$. Since $u_\lam \leq 1-\frac{3}{2}\eps$, if $\delta_\lam \leq \frac{\eps}{2}$, then $u_\lambda +\phi\leq 1-\eps$ and we have that:
\begin{equation}\arraycolsep=1.5pt\begin{array}{lll}\arraycolsep=1.5pt
J_{\eps,\lam}(u_{\lam }+\phi)-J_{\eps,\lam}(u_{\lam
})&=&\displaystyle\frac{1}{2}\displaystyle\int _{\Omega}|\nabla
\phi|^2dx+ \displaystyle\int _{\Omega}\nabla u_{\lam }\cdot\nabla \phi dx
-\lam\displaystyle\int _{\Omega} f(x) \left( \displaystyle
\frac{1}{1-u_{\lam}-\phi}-\displaystyle\frac{1}{1-u_{\lam}}\right)\\[3mm]
&\ge & \frac{\mu_{1,\lam}}{2}\int_\Omega \phi^2-
\lam \displaystyle\int
_{\Omega} f(x) \left( \displaystyle\frac{1}{1-u_{\lam}-\phi}- \displaystyle
\frac{1}{1-u_{\lam}}-\displaystyle\frac{\phi}{(1-u_{\lam})^2}-\displaystyle \frac{\phi ^2}{(1-u_{\lam })^3} \right)\,,
\end{array}\label{2:5}
\end{equation}
where we have applied $(\ref{2:4})$. Since now
$$\big| \frac{1}{1-u_{\lam}-\phi}-\frac{1}{1-u_{\lam}}-\frac{\phi}{(1-u_{\lam})^2}-\frac{\phi
^2}{(1-u_{\lam})^3}|\leq C |\phi|^3\,$$
for some $C>0$, $(\ref{2:5} )$ gives that
$$J_{\eps,\lam}(u_{\lam }+\phi)-J_{\eps,\lam}(u_{\lam })\ge \left(\frac{\mu_{1,\lam}}{2}-C\lam \parallel f \parallel_\infty \delta_\lambda \right)\int_\Omega \phi^2 >0
$$
provided $\delta_\lambda$ is small enough. This proves that $u_{\lam}$ is a local
minimum of $J_{\eps,\lam}$ in the $\mathcal{C}^1$ topology. Since (\ref{growth}) is satisfied, we can then directly apply Theorem 1 in \cite{BN} to get that $u_{\lam}$ is a local minimum of $J_{\eps,\lam}$ in $ H^1_0(\Omega)$. \qed

\medskip \noindent Since now $f \not=0$, fix some small ball $B_{2r} \subset \Omega$ of radius $2r$, $r>0$, so that $\int_{B_r} f(x) dx >0$. Take a cut-off function $\chi$ so that $\chi=1$ on $B_r$ and $\chi=0$ outside $B_{2r}$. Let $w_\eps=(1-\eps)\chi \in H_0^1(\Omega)$. We have that:
$$J_{\eps,\lam}(w_\eps)\leq \frac{(1-\eps)^2}{2}\int_\Omega |\nabla \chi|^2dx-\frac{\lam}{\eps^2}\int_{B_r} f(x) \to -\infty$$
as $\eps \to 0$, and uniformly for $\lambda$ far away from zero. Since
$$J_{\eps,\lam}(u_\lam)= \displaystyle
\frac{1}{2}\int _{\Omega}|\nabla u_\lam|^2dx -\lam \int _{\Omega}\frac{f(x)}{1-u_\lam}dx \to
\displaystyle
\frac{1}{2}\int _{\Omega}|\nabla u^*|^2dx -\lam^* \int _{\Omega}\frac{f(x)}{1-u^*}dx$$
as $\lam \to \lam^*$, we can find that for  $\eps>0$ small, the inequality
\begin{equation} \label{giu} J_{\eps,\lam}(w_\eps)<J_{\eps,\lam}(u_\lam)\end{equation}
holds for any $\lam$ close to $\lam^*$.

\medskip \noindent Fix now $\eps>0$ small enough in order that (\ref{giu}) holds for $\lam$ close to $\lam^*$, and define
$$c_{\eps,\lam} =\displaystyle\inf_{\gamma \in \Gamma} \max_{u \in \gamma}J_{\eps,\lam}(u),$$
where
$\Gamma=\{\gamma:[0,1]\to H_0^1(\Omega); \hbox{$\gamma$ continuous and $\gamma(0)=u_\lam,\:\gamma(1)=w_\eps$} \}$. We can then apply the Mountain Pass Theorem  \cite{AR} to get a solution $U_{\eps,\lam}$ of (\ref{2:1}) for $\lam$ close to $\lam^*$, provided the Palais-Smale condition holds at level $c$. We shall now prove this (PS)-condition in the following form:
\begin{lem} \label{lemPS}
Assume that $\{w_n\}\subset H^1_0(\Omega)$ satisfies
\begin{equation}
J_{\eps,\lambda_n}(w_n)\leq C\:,\qquad
J_{\eps,\lam_n}'(w_n)\to 0 \:\: \hbox{ in }H^{-1}\, \label{2:19}
\end{equation}
for $\lambda_n \to \lambda>0$. Then the sequence $(w_n)_n$ is uniformly bounded in $H^1_0(\Omega)$ and therefore admits a convergent subsequence in $H_0^1(\Omega)$.
\end{lem}
\noindent {\bf Proof:} By (\ref{2:19}) we have that:
$$ \int _{\Omega }|\nabla w_n|^2 dx-\lam_n \int _{\Omega
} f(x) g_\varepsilon(w_n)w_n dx=o(\parallel w_n \parallel_{H_0^1})$$
as $n \to +\infty$ and then,
\begin{eqnarray*}
C & \geq& \frac{1}{2} \int_\Omega |\nabla w_n|^2 dx-\lambda_n \int_\Omega f(x)G_\eps(w_n)dx\\
&=&\left(\frac{1}{2}-\frac{1}{\theta}\right)\int_\Omega |\nabla w_n|^2 dx  +\lambda_n \int_\Omega f(x)\left(\frac{1}{\theta}w_n g_\eps(w_n) -G_\eps(w_n)\right)dx+o(\parallel w_n \parallel_{H_0^1})\\
&\geq & \left(\frac{1}{2}-\frac{1}{\theta}\right)\int_\Omega |\nabla w_n|^2 dx
+\lambda_n \int_{\{w_n \geq M_\eps \}} f(x)\left(\frac{1}{\theta}w_n g_\eps(w_n) -G_\eps(w_n)\right)dx+o(\parallel w_n \parallel_{H_0^1})-C_\eps\\
&\geq & \left(\frac{1}{2}-\frac{1}{\theta}\right)\int_\Omega |\nabla w_n|^2 dx
+o(\parallel w_n \parallel_{H_0^1})-C_\eps \end{eqnarray*}
in view of (\ref{cruci}). Hence, $\displaystyle\sup_{n \in \N}\parallel w_n \parallel_{H_0^1}<+\infty$.

\medskip \noindent Since $p$ is subcritical, the compactness of the embedding $H_0^1(\Omega) \hookrightarrow L^{p+1}(\Omega)$ provides that, up to a subsequence, $w_n \to w$ weakly in $H^1_0(\Omega)$ and strongly in $L^{p+1}(\Omega)$, for some $w\in H^1_0(\Omega )$. By $(\ref{2:19})$ we get that $ \int _{\Omega }|\nabla w|^2=\lam \int _{\Omega }f(x) g_\eps(w)w$, and then, by (\ref{growth}), we deduce that
$$\int _{\Omega }|\nabla (w_n- w)|^2=\int _{\Omega }|\nabla w_n|^2-\int _{\Omega }|\nabla w|^2+o(1)=
\lambda_n \int _\Omega f(x) g_\varepsilon (w_n)w_n-\lam \int _{\Omega }f(x) g_\eps(w)w+o(1)\to 0$$
as $n \to +\infty$.\qed

\bigskip \noindent To conclude the proof of Theorem 1.2, we consider for any $\lam \in (\lam^*-\delta,\lam^*)$
the mountain pass solution $U_{\eps,\lambda}$ of (\ref{2:1}) at
energy level $c_{\eps,\lam}$, where $\delta>0$ is small enough.
Since $c_{\eps,\lambda}\leq c_{\eps,\lambda^*-\delta}$ for any $\lam
\in (\lam^*-\delta,\lam^*)$, and applying again Lemma \ref{lemPS},
we get that $\parallel U_{\eps,\lam}\parallel_{H_0^1} \leq C$, for
any $\lam$ close to $\lam^*$. Then, by (\ref{growth}) and elliptic
regularity theory, we get that $U_{\eps,\lam}$ is uniformly bounded
in $C^{2,\alpha}(\bar \Omega)$ for $\lam \uparrow \lam^*$, for
$\alpha \in (0,1)$. Hence, we can extract a sequence
$U_{\eps,\lam_n}$, $\lam_n \uparrow \lam^*$, converging in $C^2(\bar
\Omega)$ to some function $U^*$, where $U^*$ is a solution for
problem (\ref{2:1}) at $\lambda=\lambda^*$. Also $u^*$ is a solution
for (\ref{2:1}) at $\lam=\lam^*$ so that $\mu_1 \left(
-\Delta-\lambda^* f(x) g_\varepsilon'(u^*) \right)=0$. By convexity
of $g_\eps(u)$, it is classical to show that $u^*$ is the unique
solution of this equation and therefore $U^*=u^*$. Since along any
convergent sequence of $U_{\eps,\lam}$ as $\lam \uparrow \lam^*$ the
limit is always $u^*$, we get that $\lim_{\lam \uparrow
\lam^*}U_{\eps,\lam}=u^*$ in $C^2(\bar \Omega)$. Therefore, since
$u^* \leq 1-2\eps$, there exists $\delta>0$ so that for any $\lam
\in (\lam^*-\delta,\lam^*)$ $U_{\eps,\lam}\leq u^*+\eps \leq 1-\eps$
and hence, $U_{\eps,\lam}$ is a solution of $(S)_\lam$. Since the
mountain pass energy level $c_{\eps,\lam}$ satisfies
$c_{\eps,\lam}>J_{\eps,\lam}(u_\lam)$, we have that $U_{\eps,\lam}
\not =u_\lam$ and then $U_{\eps,\lam}=U_\lam$ for any $\lam \in
(\lam^*-\delta,\lam^*)$. Note that by \cite{CR},  we know that
$u_\lam$, $U_\lam$ are the only solutions of $(S)_\lam$ as $\lam
\uparrow \lam^*$.

\section{Minimal branch on the ball for power-like permittivity profiles}
Let $B$ be the unit ball. Let $(\lambda_n)_n$ be such that
$\lambda_n \to  \lam \in [0,\lam^*]$ and $u_n$ be a solution of
$(S)_{\lam_n}$ on $B$ so that $\mu_{1,n}:=\mu_{1,\lam_n}(u_n) \geq
0$. By Proposition \ref{uniqueness} $u_n$ coincides with the minimal
solution $u_{\lam_n}$ and, by some symmetrization arguments, in
\cite{GG1} it is shown that the minimal solution $u_n$ is radial and
achieves its absolute maximum only at zero.

\medskip \noindent Given a permittivity profile $f(x)$ as in (\ref{assf0}), in order to get Theorem \ref{thm0} we want to show:
\begin{equation} \label{linfinity} \sup_{n \in \N} \parallel u_n \parallel_\infty <1, \end{equation}
provided $N \geq 8$ and $\alpha>\alpha_N=\frac{3N-14-4\sqrt{6}}{4+2\sqrt{6}}$. In particular, since $u_\lam$ is non decreasing in $\lam$ and
$$\sup_{\lam \in [0,\lam^*)} \parallel u_\lam \parallel_\infty <1,$$ the extremal solution $u^*=\displaystyle\lim_{\lam \uparrow \lam^*} u_\lam$ would be a solution of $(S)_{\lam^*}$ so that $\mu_{1,\lam^*}(u^*)\leq 0$. Property (\ref{minimal*}) must hold because otherwise, by Implicit Function Theorem, we could find solutions of $(S)_\lam$ for $\lam>\lam^*$.

\bigskip \noindent In order to prove (\ref{linfinity}), let us argue by contradiction. Up to a subsequence, assume that $u_n(0)=\displaystyle \max_B u_n \to 1$ as $n \to +\infty$. Since $\lam=0$ implies $u_n \to 0$ in $C^2(\bar B)$, we can assume that $\lam_n \to \lambda>0$. Let $\eps_n:=1-u_n(0) \to 0$ as $n \to +\infty$ and introduce the following rescaled function:
\begin{equation} \label{easy} U_n(y)=\frac{1-u_n \left( \eps_n^{\frac{3}{2+\alpha}} \lam_n^{-\frac{1}{2+\alpha}} y \right) }{\eps_n}\:, \:\: y \in B_n:=B_{\eps_n^{-\frac{3}{2+\alpha}} \lam_n^{\frac{1}{2+\alpha}} }(0).\end{equation}
The function $U_n$ satisfies:
\begin{equation} \label{eqeasy} \left\{ \begin{array}{ll}
\displaystyle \Delta U_n=\frac{|y|^\alpha g \left(\eps_n^{\frac{3}{2+\alpha}} \lam_n^{-\frac{1}{2+\alpha}} y \right)
}{U_n^2}& \hbox{in }B_n,\\
\displaystyle U_n(y)\geq U_n(0)=1,  & \end{array} \right.\end{equation}
and $B_n \to \R^N$ as $n \to +\infty$. We would get a contradiction to $\mu_{1,n} \geq 0$ by proving:
\begin{prop} \label{propmax} There exists a subsequence $\{ U_n \}_n$ such that $U_n\to U $ in
$C^1_{\hbox{loc}}(\R^N)$, where $U$ is a solution of the equation:
\begin{equation} \label{limiteasy} \left\{ \begin{array}{ll} \Delta U = g(0)\displaystyle\frac{|y|^\alpha}{U^2} & \hbox{in }R^N\,,\\
U(y) \ge U(0)=1  & \hbox{in }\R^N \,.\end{array} \right. \end{equation}
Moreover, there exists $\phi_n \in C_0^\infty(B)$ such that:
\begin{eqnarray*}\int_B \left(|\nabla \phi_n|^2-\frac{2\lam_n |x|^\alpha g(x)}{(1-u_n)^3}\phi_n^2 \right) <0.\end{eqnarray*} \end{prop}
\proof Let $R>0$. For $n$ large, decompose
$U_n=U_n^1+U_n^2$, where $U_n^2$  satisfies:
$$\left\{ \begin{array}{ll} \Delta U_n^2=\Delta U_n   & \hbox{in }B_R(0) \,,\\
U^2_n =0  & \hbox{on }\partial B_R(0) \,.\end{array} \right. $$
By (\ref{eqeasy}) we get that on $B_R(0)$:
$$ 0\le \Delta U_n \leq R^\alpha \parallel g \parallel_\infty,$$
and then, standard elliptic regularity theory gives that $U^2_n$ is uniformly
bounded in $C^{1,\beta}(B_R(0))$, $\beta \in (0,1)$. Up to a subsequence, we get that $U^2_n \to U^2$ in $C^1( B_R(0))$. Since $U^1_n=U_n\ge 1$ on
$\partial B_R(0)$, by harmonicity $U^1_n\ge 1$ in $B_R(0)$ and, by Harnack inequality:
$$ \sup_{B_{R/2}(0)}U_n^1 \le C_R \inf _{B_{R/2}(0)}U_n^1\le
C_R U_n^1(0)=C_R \left(1-U_n^2(0) \right)\le C_R \left(1+\sup_{n \in \N}|U^2_n(0)| \right)<\infty\,.$$
Hence, $U^1_n$ is uniformly bounded in $C^{1,\beta}(B_{R/4}(0))$, $\beta \in (0,1)$. Up to a further subsequence,
we get that $U^1_n \to U^1$ in $C^1(B_{R/4}(0))$ and then, $U_n\to U^1+U^2$ in $C^1(B_{R/4}(0))$, for any $R>0$. By a diagonal process and up to a subsequence, we find that $U_n\to U $ in
$C^1_{\hbox{loc}}(\R^N)$, where $U$ is a solution of the equation (\ref{limiteasy}).

\medskip \noindent If $1\leq N \leq7$ or $N\geq 8$, $\alpha>\alpha_N$, since $g(0)>0$ by Theorem \ref{thm3} we get that $\mu_1(U)<0$ and then, we find $\phi \in C_0^\infty(\R^N)$ so that:
$$\int \left( |\nabla \phi|^2-2g(0)\frac{|y|^\alpha}{U^3} \phi^2\right)<0.$$
Define now $\phi_n(x)=\big( \eps_n^{\frac{3}{2+\alpha}} \lam_n^{-\frac{1}{2+\alpha}} \big)^{-\frac{N-2}{2}} \phi\left( \eps_n^{-\frac{3}{2+\alpha}} \lam_n^{\frac{1}{2+\alpha}} x\right)$. We have that:
\begin{eqnarray*}\int_B \left(|\nabla \phi_n|^2-\frac{2\lam_n |x|^\alpha g(x)}{(1-u_n)^3}\phi_n^2 \right)
&=&\int \left(|\nabla \phi|^2-\frac{2 |y|^\alpha }{U_n^3}g (\eps_n^{\frac{3}{2+\alpha}} \lam_n^{-\frac{1}{2+\alpha}}y) \phi^2 \right)\\
&\to & \int \left( |\nabla \phi|^2-2g(0)\frac{|y|^\alpha}{U^3} \phi^2\right)<0\end{eqnarray*}
as $n \to +\infty$, since $\phi$ has compact support and $U_n \to U$ in $C^1_{\hbox{loc}}(\R)$. The proof of Proposition \ref{propmax} is now complete.\qed

\section{Compactness along the second branch of solutions}
In this Section we turn  to the compactness result stated in Theorem \ref{thm2}. Assume that $f \in C(\bar \Omega)$ is in the form (\ref{assf1}), and let $(u_n)_n$ be a solution sequence for $(S)_{\lambda_n}$ where $\lambda_n \to  \lam \in[0,\lam ^*]$.

\subsection{Blow-up analysis} \label{sectionblow}
Assume that the sequence $(u_n)_n$ is not compact, which means that
up to passing to a subsequence, we may assume that $\displaystyle \max_{\Omega} u_n \to 1$ as $n\to \infty $. Let $x_n$ be a maximum point of $u_n$ in $\Omega$ (i.e., $u_n(x_n)=\displaystyle \max_{\Omega} u_n$) and set
$\varepsilon_n=1-u_n(x_n)$. Let us assume that $x_n \to p$ as $n \to +\infty$. We have three different situations depending on the location of $p$ and the rate of $|x_n-p|$:

\medskip \noindent 1) blow up outside the zero set of $f(x)$ $\{p_1,\dots,p_k\}$, i.e. $p \not \in \{p_1,\dots,p_k\}$;\\
2) ``slow'' blow up at some $p_i$ in the zero set of $f(x)$, i.e. $x_n \to p_i$ and $\eps_n^{-3}\lam_n |x_n-p_i|^{\alpha+2} \to +\infty$ as $n \to +\infty$;\\
3) ``fast'' blow at some $p_i$ in the zero set of $f(x)$, i.e. $x_n \to p_i$ and $\displaystyle\limsup_{n \to +\infty} \left( \eps_n^{-3} \lam_n |x_n-p_i|^{\alpha+2} \right)<+\infty$.

\medskip \noindent Accordingly, we discuss now each one of these situations.

\medskip \noindent {\bf $1^{st}$ Case} Assume that $p \notin \{p_1,\dots,p_k\}$. In general, we are not able to prove that a blow up point $p$ is always far away from $\partial \Omega$, even though we suspect it  to be true. However, some weaker estimate is available and --as explained later-- will be sufficient for our purposes. We have that:
\begin{lem}
Let $h_n$ be a function on  a smooth bounded domain $A_n$ in $\R^N$. Let $W_n$ be a solution of:
\begin{equation}\label{wn} \left\{ \begin{array}{ll} \displaystyle \Delta W_n=\frac{h_n(x)}{W_n^2} &\hbox{in }A_n,\\
W_n(y)\geq C>0 &\hbox{in }A_n,\\
W_n(0)=1, & \end{array}\right.\end{equation}
for some $C>0$. Assume that $\displaystyle\sup_{n \in \N}\parallel h_n \parallel_\infty <+\infty$ and $A_n \to T_\mu$ as $n \to +\infty$ for some $\mu\in (0,+\infty)$, where $T_\mu$ is an hyperspace so that $0 \in T_\mu$ and $\hbox{dist }(0,\partial T_\mu)=\mu$. Then, either
\begin{equation}\label{condition1} \displaystyle \inf_{\partial A_n \cap B_{2\mu}(0)} W_n \leq C \end{equation}
or
\begin{equation}\label{condition2} \inf_{\partial A_n \cap B_{2\mu}(0)} \partial_{\nu} W_n \leq 0, \end{equation}
where $\nu$ is the unit outward normal of $A_n$. \label{bconc} \end{lem}
\proof Assume that $\partial_{\nu} W_n>0$ on $\partial A_n \cap B_{2\mu}(0)$. Let
$$G(x)=\left\{ \begin{array}{ll }-\frac{1}{2\pi}\log \frac{|x|}{2\mu}& \hbox{if }N=2\\
c_N \left(\frac{1}{|x|^{N-2}}-\frac{1}{(2\mu)^{N-2}}\right) &\hbox{if }N\geq 3 \end{array}\right.$$
be the Green function at $0$ of the operator $-\Delta$ in $B_{2\mu}(0)$ with homogeneous Dirichlet boundary condition, where $c_N=\frac{1}{(N-2)|\partial B_1(0)|}$ and $|\cdot|$ stands for the Lebesgue measure.

\medskip \noindent Here and in the sequel, when there is no ambiguity on the domain we are considering, $\nu$ and $d\sigma$ will denote the unit outward normal and the boundary integration element of the corresponding domain. By the representation formula we have that:
\begin{eqnarray}
W_n(0)&=&-\int_{A_n \cap B_{2\mu}(0)}\Delta W_n(x) G(x) dx -\int_{\partial A_n \cap B_{2\mu}(0)}W_n(x) \partial_{\nu}G(x) d\sigma(x) \label{repr}\\
&&+\int_{\partial A_n \cap B_{2\mu}(0)} \partial_{\nu} W_n(x) G(x) d\sigma(x)-\int_{\partial B_{2\mu}(0) \cap A_n} W_n(x) \partial_{\nu}G(x) d \sigma(x). \nonumber \end{eqnarray}
Since on $\partial T_\mu$:
\begin{equation} \label{derivbordo}-\partial_{\nu}G(x)=\left\{ \begin{array}{ll} \frac{1}{2\pi}\frac{x}{|x|^2}\cdot \nu & \hbox{if }N=2\\
(N-2) c_N \frac{x}{|x|^N}\cdot \nu &\hbox{if }N\geq 3   \end{array} \right. >0 \end{equation}
and $\partial A_n \to \partial T_\mu$, we get that
\begin{equation} \label{derivbordo1} \partial_{\nu}G(x)<0  \quad \hbox{on } \partial A_n\cap B_{2 \mu}(0).\end{equation}
Hence, by (\ref{repr}), (\ref{derivbordo1}) and the assumptions on $W_n$, we then get:
\begin{eqnarray*} 1 &\geq& -\int_{A_n \cap B_{2\mu}(0)} \frac{h_n(x)}{W_n^2(x)} G(x) dx
- \left(\inf_{\partial A_n \cap B_{2\mu}(0)} W_n \right)  \int_{\partial A_n \cap B_{2\mu}(0)} \partial_{\nu}G(x) d\sigma(x)
\end{eqnarray*}
since $G(x)\geq 0$ in $B_{2\mu}(0)$ and $\partial_{\nu}G(x)\leq 0$ on $\partial B_{2\mu}(0)$. Now, we have that
$$\big| \int_{A_n \cap B_{2\mu}(0)} \frac{h_n(x)}{W_n^2(x)} G(x) \big|\leq C$$
and by (\ref{derivbordo})
$$-\int_{\partial A_n \cap B_{2\mu}(0)} \partial_{\nu}G(x) d\sigma(x) \to -\int_{\partial T_\mu \cap B_{2\mu}(0)} \partial_{\nu}G(x) d\sigma(x)>0.$$
Then, $1\geq -C+C^{-1} \left(\displaystyle \inf_{\partial A_n \cap B_{2\mu}(0)} W_n \right)$ for some $C>0$ large enough. Hence, $\displaystyle \inf_{\partial A_n \cap B_{2\mu}(0)} W_n$ is uniformly bounded and the proof is complete. \qed

\medskip \noindent We are now ready to completely discuss this first case. Introduce the following rescaled function:
\begin{equation} \label{un1caso} U_n(y)=\frac{1-u_n(\varepsilon_n^{\frac{3}{2}}\lam_n^{-\frac{1}{2}}y+x_n)}{\varepsilon _n}\:, \qquad  y \in \Omega _n=\frac{\Omega
-x_n}{\varepsilon_n^{ \frac{3}{2} }\lam_n^{-\frac{1}{2}}} \,.\end{equation}
Then, $U_n$ satisfies
\begin{equation} \label{gente1caso} \left\{ \begin{array}{ll}
\displaystyle \Delta U_n =\frac{f(\varepsilon _n^{\frac{3}{2}}\lam_n^{-\frac{1}{2}}y+x_n)}{U_n^2} & \hbox{in }\Omega _n,\\
\displaystyle U_n(0)=1  & \hbox{in }\Omega _n.\end{array} \right. \end{equation}
In addition, we have that $U_n \geq U_n(0)=1$ as long as $x_n$ is the maximum point of $u_n$ in $\Omega$.

\medskip \noindent We would like to prove the following:
\begin{prop} \label{prop1caso} Let $x_n \in \Omega$ and set $\eps_n:=1-u_n(x_n)$. Assume that
\begin{equation} \label{scemotta} x_n \to p \notin \{p_1,\dots,p_k \} \:, \:\: \eps_n^3 \lam_n^{-1} \to 0 \quad \hbox{as }n \to +\infty.\end{equation}
Let $U_n$, $\Omega_n$ be defined as in (\ref{un1caso}). Assume that
\begin{equation} \label{lontanozero} U_n \geq C>0 \quad \hbox{in }\Omega_n \cap B_{R_n}(0)\end{equation}
for some $R_n \to +\infty$ as $n \to +\infty$. Then, there exists a subsequence of $(U_n)_n$ such that $U_n\to U $ in
$C^1_{\hbox{loc}}(\R^N)$, where $U$ is a solution of the equation:
\begin{equation} \label{limitun1caso} \left\{ \begin{array}{ll} \Delta U = \displaystyle\frac{f(p)}{U^2} & \hbox{in }R^N\,,\\
U(y) \geq C>0  & \hbox{in }\R^N \,.\end{array} \right. \end{equation}
Moreover, there exists a function $\phi_n \in C_0^\infty(\Omega)$ such that:
\begin{equation} \label{testfunction} \int_\Omega \left(|\nabla \phi_n|^2-\frac{2\lam_n f(x)}{(1-u_n)^3}\phi_n^2 \right)<0 \end{equation}
and $\hbox{Supp }\phi_n \subset B_{M \eps_n^{\frac{3}{2}} \lam_n^{-\frac{1}{2}}}(x_n)$ for some $M>0$.
\end{prop}
\proof By (\ref{scemotta}) Lemma \ref{bconc} provides us with a stronger estimate:
\begin{equation} \eps_n^3 \lam_n^{-1} (\hbox{dist }(x_n,\partial \Omega))^{-2} \to 0 \quad \hbox{as }n \to +\infty.\label{speed} \end{equation}
Indeed, by contradiction and up to a subsequence, assume that $\eps_n^3 \lam_n^{-1} d_n^{-2} \to \delta>0$ as $n \to +\infty$, where $d_n:=\hbox{dist }(x_n,\partial \Omega)$. In view of (\ref{scemotta}) we get that $d_n \to 0$ as $n \to +\infty$. We introduce the following rescaling $W_n$:
$$W_n(y)=\frac{1-u_n(d_n y+x_n)}{\eps_n}\,,\quad y \in A_n=\frac{\Omega-x_n}{d_n} \,. $$
Since $d_n \to 0$, we get that $A_n \to T_1$ as $n \to +\infty$, where $T_\mu$ is an hyperspace containing $0$ so that $\hbox{dist }(0,\partial T_\mu)=\mu$. The function $W_n$ solves problem (\ref{wn}) with $h_n(y)=\frac{\lam_n d_n^2}{\eps_n^3}f(d_n y+x_n)$ and $C=W_n(0)=1$.
We have that:
$$\parallel h_n \parallel_\infty \leq \frac{\lam_n d_n^2}{\eps_n^3} \parallel f \parallel_\infty \leq \frac{2}{\delta}\parallel f \parallel_\infty$$
and $W_n=\frac{1}{\eps_n} \to + \infty$ on $\partial A_n$. By Lemma \ref{bconc} we get that (\ref{condition2}) must hold. A contradiction to Hopf Lemma applied to $u_n$. Hence, the validity of (\ref{speed}).

\medskip \noindent We have proved that the blow up is ``essentially'' in the interior of $\Omega$: (\ref{speed}) implies that $\Omega_n \to \R^N$ as $n \to +\infty$. Arguing as in the proof of Proposition \ref{propmax}, we get that $U_n \to U$ in $C^1_{\hbox{loc}}(\R^N)$, where $U$ is a solution of (\ref{limitun1caso}) by means of (\ref{gente1caso})-(\ref{lontanozero}).

\medskip \noindent If $1\leq N \leq7$, since $f(p)>0$ by Theorem \ref{thm3} we get that $\mu_1(U)<0$ and then, we find $\phi \in C_0^\infty(\R^N)$ so that:
$$\int \left( |\nabla \phi|^2-\frac{2f(p)}{U^3} \phi^2\right)<0.$$
Define now $\phi_n(x)=\big( \eps_n^{\frac{3}{2}} \lam_n^{-\frac{1}{2}} \big)^{-\frac{N-2}{2}} \phi\left( \eps_n^{-\frac{3}{2}} \lam_n^{\frac{1}{2}} (x-x_n)\right)$. We have that:
\begin{eqnarray*}\int_\Omega \left(|\nabla \phi_n|^2-\frac{2\lam_n f(x)}{(1-u_n)^3}\phi_n^2 \right)
&=&\int \left(|\nabla \phi|^2-\frac{2  f (\eps_n^{\frac{3}{2}} \lam_n^{-\frac{1}{2}}y+x_n) }{U_n^3} \phi^2 \right)\\
&\to & \int \left( |\nabla \phi|^2-\frac{2 f(p)}{U^3} \phi^2\right)<0\end{eqnarray*}
as $n \to +\infty$, since $\phi$ has compact support and $U_n \to U$ in $C^1_{\hbox{loc}}(\R)$. The proof of Proposition \ref{prop1caso} is now complete.\qed

\medskip \noindent {\bf $2^{nd}$ Case} Assume that $x_n \to p_i$ and $\eps_n^{-3}\lam_n |x_n-p_i|^{\alpha+2} \to +\infty$ as $n \to +\infty$. Define
\begin{equation} \label{deffi}  f_i(x)=\left(\prod_{j=1,\:j \not=i}^k |x-p_j|^{\alpha_j}\right) g(x).\end{equation}
We rescale the function $u_n$ in a different way:
\begin{equation} \label{un2caso} U_n(y)=\frac{1-u_n(\eps_n ^{\frac{3}{2}}\lam_n^{-\frac{1}{2}}|x_n-p_i|^{-\frac{\alpha}{2}} y+x_n)}{\varepsilon _n}\:, \qquad  y \in \Omega _n=\frac{\Omega-x_n}{\eps_n ^{\frac{3}{2}}\lam _n^{-\frac{1}{2}}|x_n-p_i|^{-\frac{\alpha}{2}}} \,.\end{equation}
In this situation, $U_n$ satisfies:
\begin{equation} \label{gente2caso} \left\{ \begin{array}{ll}
\displaystyle \Delta U_n =\big| \eps_n^{\frac{3}{2}} \lam_n^{-\frac{1}{2}} |x_n-p_i|^{-\frac{\alpha+2}{2}}y +\frac{x_n-p_i}{|x_n-p_i|} \big|^\alpha \frac{f_i(\eps_n^{\frac{3}{2}}\lam_n^{-\frac{1}{2}}|x_n-p_i|^{-\frac{\alpha}{2}}y+x_n)}{U_n^2} & \hbox{in }\Omega _n,\\
\displaystyle U_n(0)=1  & \hbox{in }\Omega _n.\end{array} \right. \end{equation}

\medskip \noindent The following result holds:
\begin{prop} \label{prop2caso} Let $x_n \in \Omega$ and set $\eps_n:=1-u_n(x_n)$. Assume that
\begin{equation} \label{scemotta2caso} x_n \to p_i \:,\:\: \eps_n^{-3}\lam_n |x_n-p_i|^{\alpha+2} \to +\infty \quad \hbox{as }n \to +\infty.\end{equation}
Let $U_n$, $\Omega_n$ be defined as in (\ref{un2caso}). Assume that (\ref{lontanozero}) holds. Then, up to a subsequence, $U_n\to U $ in $C^1_{\hbox{loc}}(\R^N)$, where $U$ is a solution of the equation:
\begin{equation} \label{limitun2caso} \left\{ \begin{array}{ll} \Delta U = \displaystyle\frac{f_i(p_i)}{U^2} & \hbox{in }R^N\,,\\
U(y) \geq C>0  & \hbox{in }\R^N \,.\end{array} \right. \end{equation}
Moreover, there holds (\ref{testfunction}) for some $\phi_n \in C_0^\infty(\Omega)$ such that $\hbox{Supp }\phi_n \subset B_{M \eps_n^{\frac{3}{2}} \lam_n^{-\frac{1}{2}} |x_n-p_i|^{-\frac{\alpha}{2}} }(x_n)$, $M>0$.
\end{prop}
\proof By (\ref{scemotta2caso}) we get that $\Omega_n \to \R^N$ as $n \to +\infty$. As before, $U_n \to U$ in $C^1_{\hbox{loc}}(\R^N)$ and $U$ is a solution of (\ref{limitun2caso}) in view of (\ref{lontanozero}) and (\ref{gente2caso})-(\ref{scemotta2caso}). Since $1\leq N \leq7$ and $f_i(p_i)>0$, Theorem \ref{thm3} implies $\mu_1(U)<0$ and the existence of some $\phi \in C_0^\infty(\R^N)$ so that:
$$\int \left( |\nabla \phi|^2-\frac{2 f_i (p_i)}{U^3} \phi^2\right)<0.$$
Define now $\phi_n(x)=\big( \eps_n^{\frac{3}{2}} \lam_n^{-\frac{1}{2}}|x_n-p_i|^{-\frac{\alpha}{2}} \big)^{-\frac{N-2}{2}} \phi\left( \eps_n^{-\frac{3}{2}} \lam_n^{\frac{1}{2}} |x_n-p_i|^{\frac{\alpha}{2}} (x-x_n)\right)$. We have that:
\begin{eqnarray*}&&\int_\Omega \left(|\nabla \phi_n|^2-\frac{2\lam_n f(x)}{(1-u_n)^3}\phi_n^2 \right)\\
&& =\int \left(|\nabla \phi|^2-\big| \eps_n^{\frac{3}{2}} \lam_n^{-\frac{1}{2}} |x_n-p_i|^{-\frac{\alpha+2}{2}}y +\frac{x_n-p_i}{|x_n-p_i|} \big|^\alpha
\frac{2  f_i (\eps_n^{\frac{3}{2}} \lam_n^{-\frac{1}{2}} |x_n-p_i|^{-\frac{\alpha}{2}} y+x_n) }{U_n^3} \phi^2 \right)\\
&&\to  \int \left( |\nabla \phi|^2-\frac{2 f_i(p_i)}{U^3} \phi^2\right)<0\end{eqnarray*}
as $n \to +\infty$. Proposition \ref{prop2caso} is now completely proved.\qed

\medskip \noindent {\bf $3^{rd}$ Case} Assume that $x_n \to p_i$ as $n \to +\infty$ and $\eps_n^{-3}\lam_n |x_n-p_i|^{\alpha+2} \leq C$. We rescale the function $u_n$ in a still different way:
\begin{equation} \label{un3caso} U_n(y)=\frac{1-u_n(\eps_n^{\frac{3}{2+\alpha}}\lam_n^{-\frac{1}{2+\alpha}} y+x_n)}{\varepsilon _n}\:, \qquad  y \in \Omega_n=\frac{\Omega-x_n}{\eps_n^{\frac{3}{2+\alpha}}\lam_n^{-\frac{1}{2+\alpha}}} \,.\end{equation}
The equation satisfied by $U_n$ is:
\begin{equation} \label{gente3caso} \left\{ \begin{array}{ll}
\displaystyle \Delta U_n =\big|y+ \eps_n^{-\frac{3}{2+\alpha}} \lam_n^{\frac{1}{2+\alpha}} (x_n-p_i) \big|^\alpha
\frac{f_i(\eps_n^{\frac{3}{2+\alpha}}\lam_n^{-\frac{1}{2+\alpha}}y+x_n)}{U_n^2} & \hbox{in }\Omega _n,\\
\displaystyle U_n(0)=1  & \hbox{in }\Omega _n,\end{array} \right. \end{equation}
where $f_i$ is defined in (\ref{deffi}).

\medskip \noindent In this situation, the result we have is the following:
\begin{prop} \label{prop3caso} Let $x_n \in \Omega$ and set $\eps_n:=1-u_n(x_n)$. Assume that
\begin{equation} \label{scemotta3caso} \eps_n^3 \lam_n^{-1} \to 0 \:,\:\: x_n \to p_i \:,\:\: \eps_n^{-\frac{3}{\alpha+2}} \lam_n^{\frac{1}{\alpha+2}} (x_n-p_i) \to y_0 \quad \hbox{as }n \to +\infty.\end{equation}
Let $U_n$, $\Omega_n$ be defined as in (\ref{un3caso}). Assume that either (\ref{lontanozero}) holds or
\begin{equation} \label{lontanozerobis}
U_n \geq C \left(\eps_n^{-\frac{3}{\alpha+2}}  \lam_n^{\frac{1}{\alpha+2}} |x_n-p_i| \right)^{-\frac{\alpha}{3}}
\big| y+\eps_n^{-\frac{3}{\alpha+2}}\lam_n^{\frac{1}{\alpha+2}}(x_n-p_i) \big|^{\frac{\alpha}{3}}
\quad \hbox{in }\Omega_n \cap B_{R_n}(0)\end{equation}
for some $R_n \to +\infty$ as $n \to +\infty$ and $C>0$. Then, up to a subsequence, $U_n\to U $ in $C^1_{\hbox{loc}}(\R^N)$ and $U$ satisfies:
\begin{equation} \label{limitun3caso} \left\{ \begin{array}{ll} \Delta U = |y+y_0|^\alpha \displaystyle\frac{f_i(p_i)}{U^2} & \hbox{in }\R^N\,,\\
U(y) \geq C>0  & \hbox{in }\R^N \,.\end{array} \right. \end{equation}
Moreover, we have that (\ref{testfunction}) holds for some function $\phi_n \in C_0^\infty(\Omega)$ such that $\hbox{Supp }\phi_n \subset B_{M \eps_n^{\frac{3}{2+\alpha}} \lam_n^{-\frac{1}{2+\alpha}}}(x_n)$, $M>0$.
\end{prop}
\proof By (\ref{scemotta3caso}) we get that $\Omega_n \to \R^N$ as $n \to +\infty$. If (\ref{lontanozero}) holds, as before $U_n \to U$ in $C^1_{\hbox{loc}}(\R^N)$ and, by (\ref{lontanozero}) and (\ref{gente3caso})-(\ref{scemotta3caso}), $U$ solves (\ref{limitun3caso}).

\medskip \noindent We need to discuss the non trivial situation when we have the validity of (\ref{lontanozerobis}). Arguing as in the proof of Proposition \ref{propmax}, fix $R>2|y_0|$ and decompose $U_n=U_n^1+U_n^2$, where $U_n^2$  satisfies:
$$\left\{ \begin{array}{ll} \Delta U_n^2=\Delta U_n   & \hbox{in }B_R(0) \,,\\
U^2_n =0  & \hbox{on }\partial B_R(0) \,.\end{array} \right. $$
By (\ref{gente3caso}) and (\ref{lontanozerobis}) we get that on $B_R(0)$:
\begin{eqnarray*} 0 \le \Delta U_n &=& \big|y+ \eps_n^{-\frac{3}{2+\alpha}} \lam_n^{\frac{1}{2+\alpha}} (x_n-p_i) \big|^\alpha
\frac{f_i(\eps_n^{\frac{3}{2+\alpha}}\lam_n^{-\frac{1}{2+\alpha}}y+x_n)}{U_n^2}\\
&\leq& C \left(\eps_n^{-\frac{3}{\alpha+2}}  \lam_n^{\frac{1}{\alpha+2}} |x_n-p_i| \right)^{\frac{2\alpha}{3}} \big|y+ \eps_n^{-\frac{3}{2+\alpha}} \lam_n^{\frac{1}{2+\alpha}} (x_n-p_i) \big|^{\frac{\alpha}{3}}.
\end{eqnarray*}
Since $\eps_n^{-\frac{3}{\alpha+2}} \lam_n^{\frac{1}{\alpha+2}} (x_n-p_i)$ is bounded, we get that $0\leq \Delta U_n \leq C_R$ on $B_R(0)$ for $n$ large, and then, standard elliptic regularity theory gives that $U^2_n$ is uniformly
bounded in $C^{1,\beta}(B_R(0))$, $\beta \in (0,1)$. Up to a subsequence, we get that $U^2_n \to U^2$ in $C^1( B_R(0))$. Since by (\ref{lontanozerobis}) $U^1_n=U_n\ge C(R-2|y_0|)^{\frac{\alpha}{3}}>0$ on
$\partial B_R(0)$, by harmonicity $U^1_n\ge C_R$ in $B_R(0)$ and, by Harnack inequality:
$$ \sup_{B_{R/2}(0)}U_n^1 \le C_R \inf _{B_{R/2}(0)}U_n^1\le
C_R U_n^1(0)=C_R \left(1-U_n^2(0) \right)\le C_R \left(1+\sup_{n \in \N}|U^2_n(0)| \right)<\infty\,.$$
Hence, $U^1_n$ is uniformly bounded in $C^{1,\beta}(B_{R/4}(0))$, $\beta \in (0,1)$. Up to a further subsequence,
we get that $U^1_n \to U^1$ in $C^1(B_{R/4}(0))$ and then, $U_n\to U^1+U^2$ in $C^1(B_{R/4}(0))$, for any $R>0$. By a diagonal process and up to a subsequence, by (\ref{lontanozerobis}) we find that $U_n\to U $ in
$C^1_{\hbox{loc}}(\R^N)$, where $U \in C^1(\R^N) \cap C^2(\R^N \setminus \{-y_0 \})$ is a solution of the equation
$$ \left\{ \begin{array}{ll} \Delta U = |y+y_0|^\alpha \displaystyle\frac{f_i(p_i)}{U^2} & \hbox{in }\R^N \setminus \{-y_0 \}  \,,\\
U(y) \geq C|y+y_0|^{\frac{\alpha}{3}} & \hbox{in }\R^N \,,\end{array} \right.$$
for some $C>0$. In order to prove that $U$ is a solution of (\ref{limitun3caso}), we need to prove that $U(-y_0)>0$. Let $B$ some ball so that $-y_0 \in \partial B$ and assume by contradiction that $U(-y_0)=0$. Since
$$-\Delta U+c(y) U=0 \hbox{ in }B\:,\:\: U \in C^2(B)\cap C(\bar B)\:,\:\: U(y)>U(-y_0) \hbox{ in }B,$$
and $c(y)=f_i(p_i) \frac{|y+y_0|^\alpha}{U^3} \geq 0$ is a bounded function, by Hopf Lemma we get that $\partial_\nu U(-y_0)<0$, where $\nu$ is the unit outward normal of $B$ at $-y_0$. Hence, $U$ becomes negative in a neighborhood of $-y_0$ in contradiction with the positivity of $U$. Hence, $U(-y_0)>0$ and $U$ satisfies (\ref{limitun3caso}).

\medskip \noindent Since $1\leq N \leq7$ and $f_i(p_i)>0$, Theorem \ref{thm3} implies $\mu_1(U)<0$ and the existence of some $\phi \in C_0^\infty(\R^N)$ so that:
$$\int \left( |\nabla \phi|^2- |y+y_0|^\alpha \frac{2 f_i(p_i)}{U^3} \phi^2\right)<0.$$
Let $\phi_n(x)=\big( \eps_n^{\frac{3}{2+\alpha}} \lam_n^{-\frac{1}{2+\alpha}} \big)^{-\frac{N-2}{2}} \phi\left( \eps_n^{-\frac{3}{2+\alpha}} \lam_n^{\frac{1}{2+\alpha}} (x-x_n)\right)$. There holds:
\begin{eqnarray*}\int_\Omega \left(|\nabla \phi_n|^2-\frac{2\lam_n f(x)}{(1-u_n)^3}\phi_n^2 \right)
&=&\int \left(|\nabla \phi|^2-\big| y+ \eps_n^{-\frac{3}{2+\alpha}} \lam_n^{\frac{1}{2+\alpha}} (x_n-p_i) \big|^\alpha
\frac{2  f_i (\eps_n^{\frac{3}{2+\alpha}} \lam_n^{-\frac{1}{2+\alpha}}y+x_n) }{U_n^3} \phi^2 \right)\\
&\to & \int \left( |\nabla \phi|^2-|y+y_0|^\alpha \frac{2 f_i(p_i)}{U^3} \phi^2\right)<0\end{eqnarray*}
as $n \to +\infty$. Also Proposition \ref{prop3caso} is established.\qed

\subsection{Spectral confinement}
Let us assume now the validity of (\ref{morse}), namely $\mu_{2,n}:=\mu_{2,\lambda_n}(u_n)\geq 0$ for any $n \in \N$. This information will play a crucial role in  controlling the number $k$ of ``blow up points'' (for $(1-u_n)^{-1}$) in terms of the spectral information on $u_n$. Indeed, roughly speaking, we can estimate $k$ with the number of negative eigenvalues of $L_{u_n,\lam_n}$ (with multiplicities). In particular, assumption (\ref{morse}) implies that ``blow up'' can occur only along the sequence $x_n$ of maximum points of $u_n$ in $\Omega$. The following pointwise estimate on $u_n$ is available:
\begin{prop} \label{propest} Assume $2\leq N \leq 7$. Let $f \in C(\bar \Omega)$ be as in (\ref{assf1}). Let $\lam_n \to \lam \in [0,\lam^*]$ and $u_n$ be an associated solution. Assume that $u_n(x_n)=\displaystyle \max_{\Omega}u_n \to 1$ as $n \to +\infty$. Then, there exist constants $C>0$ and $N_0\in  \N$ such that
\begin{equation}
\big( 1-u_n(x)\big)\ge C \lam _n^{\frac{1}{3}} d(x)^{\frac{\alpha}{3}} |x-x_n|^{\frac{2}{3}}
\,,\quad \forall \: x \in \Omega\:,\:\: n\ge N_0\,,  \label{3:16}
\end{equation}
where $d(x)=\min \{|x-p_i|:\:i=1,\dots,k \}$ is the distance function from the zero set of $f(x)$ $\{p_1,\dots,p_k\}$.
\end{prop}
\proof Let $\eps_n=1-u_n(x_n)$. Then, $\eps_n \to 0$ as $n \to +\infty$ and, even more precisely:
\begin{equation}\label{behave} \eps_n^2 \lam_n^{-1} \to 0 \quad \hbox{as } n \to +\infty.\end{equation}
Indeed otherwise, we would have along some subsequence:
$$0 \leq \displaystyle\frac{\lam_n f(x)}{(1-u_n)^2}\leq \displaystyle\frac{\lam_n}{\eps_n^2}\parallel f \parallel_\infty \leq C\:,\:\: \lam_n \to 0 \hbox{ as }n \to +\infty.$$
But if the right hand side of $(S)_{\lam_n}$ is uniformly bounded, from elliptic regularity theory we get that $u_n$ is uniformly bounded in $C^{1,\beta}(\bar \Omega)$, $\beta \in (0,1)$. Hence, up to a furhter subsequence, $u_n \to u$ in $C^1(\bar \Omega)$, where $u$ is an harmonic function such that $u=0$ on $\partial \Omega$, $\displaystyle\max_{\Omega}u=1$. A contradiction.

\medskip \noindent By (\ref{behave}) we get that $\eps_n^3 \lam_n^{-1} \to 0$ as $n \to +\infty$, as needed in (\ref{scemotta}), (\ref{scemotta3caso}) respectively in Proposition \ref{prop1caso}, \ref{prop3caso}. Now, depending on the case corresponding to the blow up sequence $x_n$, we can apply one among Propositions \ref{prop1caso}-\ref{prop3caso} to get the existence of a function $\phi_n \in C_0^\infty(\Omega)$ such that (\ref{testfunction}) holds and with a specific control on $\hbox{Supp }\phi_n$.

\medskip \noindent By contradiction, assume now that (\ref{3:16}) is false: up to a subsequence, there exist a sequence $y_n \in \Omega$  such that
\begin{equation}
\lam_n^{-\frac{1}{3}} d(y_n)^{-\frac{\alpha}{3}} |y_n-x_n|^{-\frac{2}{3}}\big(
1-u_n(y_n)\big)=\lam_n^{-\frac{1}{3}}\min_{x \in \Omega
}\Big( d(x)^{-\frac{\alpha}{3}} |x-x_n|^{-\frac{2}{3}}\left( 1-u_n(x) \right) \Big)  \to 0\quad as
\quad n\to \infty\,. \label{3:18}
\end{equation}
Then, $\mu_n:=1-u_n(y_n) \to 0$ as $ n\to \infty$ and (\ref{3:18}) rewrites as:
\begin{equation} \label{conseq} \frac{\mu_n^{\frac{3}{2}} \lam_n^{-\frac{1}{2}}}{|x_n-y_n| \: d(y_n)^{\frac{\alpha}{2}} } \to 0 \quad \hbox{as }  n \to +\infty.
\end{equation}
We want now to explain the meaning of the crucial choice (\ref{3:18}). Let $\beta_n$ be a sequence of positive numbers so that
\begin{equation} \label{assRn}R_n:=\beta_n^{-\frac{1}{2}}\min \{d(y_n)^{\frac{1}{2}},|x_n-y_n|^{\frac{1}{2}} \} \to +\infty \quad \hbox{as } n \to +\infty.\end{equation}
Let us introduce the following rescaled function:
$$\hat U_n(y)=\frac{1-u_n(\beta_n y+y_n)}{\mu_n}\:,\:\:  y \in \hat \Omega_n=\frac{\Omega-y_n}{\beta_n}.$$
Formula (\ref{3:18}) implies:
\begin{eqnarray*}  \mu_n &=& d(y_n)^{\frac{\alpha}{3}} |y_n-x_n|^{\frac{2}{3}} \min_{x \in \Omega
}\Big( d(x)^{-\frac{\alpha}{3}} |x-x_n|^{-\frac{2}{3}}\left( 1-u_n(x) \right) \Big)\\
&\leq&  \mu_n  d(y_n)^{\frac{\alpha}{3}} |y_n-x_n|^{\frac{2}{3}} d(\beta_n y+y_n)^{-\frac{\alpha}{3}} |\beta_n y+y_n-x_n|^{-\frac{2}{3}} \hat U_n(y).\end{eqnarray*}
Since
$$\frac{d(\beta_n y+y_n)}{d(y_n)}=\min \{\big| \frac{y_n-p_i}{d(y_n)}+\frac{\beta_n}{d(y_n)}y \big|:\:i=1,\dots,k  \}\geq 1-\frac{\beta_n}{d(y_n)}|y|$$
in view of $|y_n-p_i|\geq d(y_n)$, by (\ref{assRn}) we get that:
$$\hat U_n (y) \ge  \left(1- \displaystyle \frac{\beta_n R_n}{d(y_n)} \right)^{\frac{\alpha}{3}}
\left(1- \displaystyle \frac{\beta_n R_n}{|x_n-y_n|} \right)^{\frac{2}{3}}
\ge \big(\displaystyle\frac{1}{2}\big)^{\frac{2+\alpha}{3}}$$
for any $y \in \hat \Omega_n \cap B_{R_n}(0)$. Hence, whenever (\ref{assRn}) holds, we get the validity of (\ref{lontanozero}) for the rescaled function $\hat U_n$ at $y_n$ with respect to $\beta_n$.

\medskip \noindent We need to discuss all the possible types of blow up at $y_n$.

\medskip \noindent {\bf $1^{\hbox{st}}$ Case} Assume that $y_n \to q \notin \{p_1,\dots,p_k\}$. By (\ref{conseq}) we get that $\mu_n^3 \lam_n^{-1} \to 0$ as $n \to +\infty$. Since $d(y_n)\geq C>0$, let $\beta_n=\mu_n^{\frac{3}{2}} \lam_n^{-\frac{1}{2}}$ and, by (\ref{conseq}) we get that (\ref{assRn}) holds. Associated to $y_n$, $\mu_n$, define $\hat U_n$, $\hat \Omega_n$ as in (\ref{un1caso}). We have that (\ref{lontanozero}) holds by the validity of (\ref{assRn}) for our choice of $\beta_n$. Hence, Proposition \ref{prop1caso} applies to $\hat U_n$ and give the existence of $\psi_n \in C_0^\infty(\Omega)$ such that (\ref{testfunction}) holds and $\hbox{Supp }\psi_n \subset B_{M \mu_n^{\frac{3}{2}} \lam_n^{-\frac{1}{2}}}(y_n)$, $M>0$. In the worst case $x_n \to q$, given $ U_n$ be as in (\ref{un1caso}) associated to $x_n$, $\eps_n$, we get by scaling that for $x=\eps_n^{\frac{3}{2}} \lam_n^{-\frac{1}{2}}y+x_n$:
\begin{eqnarray*}\lam_n^{-\frac{1}{3}} \big( d(x)^{-\frac{\alpha}{3}} |x-x_n|^{-\frac{2}{3}}\left( 1-u_n(x) \right)\big)\geq C
\lam_n^{-\frac{1}{3}} \big(|x-x_n|^{-\frac{2}{3}}\left( 1-u_n(x) \right)\big)= C |y|^{-\frac{2}{3}} U_n(y)\geq C_R>0\end{eqnarray*}
uniformly in $n$ and $y \in B_R(0)$, for any $R>0$. Then,
$$ \frac{\eps_n^{\frac{3}{2}} \lam_n^{-\frac{1}{2}} }{|x_n-y_n|} \to 0 \quad \hbox{as }n \to +\infty.$$
Hence, in this situation $\phi_n$ and $\psi_n$ have disjoint compact supports and obviously, it remains true when $x_n \to p \not=q$. Hence, $\mu_{2,n}<0$ in contradiction with (\ref{morse}).

\medskip \noindent {\bf $2^{\hbox{nd}}$ Case} Assume that $y_n \to p_i$ in a ``slow'' way:
$$\mu_n^{-3}\lam_n |y_n-p_i|^{\alpha+2} \to +\infty \quad \hbox {as }n \to +\infty.$$ Let now $\beta_n=\mu_n^{\frac{3}{2}} \lam_n^{-\frac{1}{2}} |y_n-p_i|^{-\frac{\alpha}{2}}$. Since $d(y_n)=|y_n-p_i|$ in this situation, we get that:
$$\frac{d(y_n)}{\beta_n}= \mu_n^{-\frac{3}{2}} \lam_n^{\frac{1}{2}}|y_n-p_i|^{\frac{\alpha+2}{2}} \to +\infty,$$
and (\ref{conseq}) gives exactly:
\begin{equation} \label{disgiunte} \frac{|x_n-y_n|}{\beta_n}=\frac{|x_n-y_n|}{\mu_n^{\frac{3}{2}} \lam_n^{-\frac{1}{2}} |y_n-p_i|^{-\frac{\alpha}{2}}} \to +\infty \end{equation}
as $n \to +\infty$. Hence, (\ref{assRn}) holds. Associated to $\mu_n$, $y_n$, define now $\hat U_n$, $\hat \Omega_n$ according to (\ref{un2caso}). Since (\ref{lontanozero}) follows by (\ref{assRn}), Proposition \ref{prop2caso} for $\hat U_n$ gives  some $\psi_n \in C_0^\infty(\Omega)$ such that (\ref{testfunction}) holds and $\hbox{Supp }\psi_n \subset B_{M \mu_n^{\frac{3}{2}} \lam_n^{-\frac{1}{2}} |y_n-p_i|^{-\frac{\alpha}{2}} }(y_n)$, $M>0$. If $x_n \to p \not=p_i$, then clearly $\phi_n$, $\psi_n$ have disjoint compact supports leading to $\mu_{2,n}<0$ in contradiction with (\ref{morse}). If also $x_n \to p_i$, we can easily show by scaling that:\\
1) if $\eps_n^{-3}\lam_n |x_n-p_i|^{\alpha+2} \to +\infty$ as $n \to +\infty$, given $ U_n$ be as in (\ref{un2caso}) associated to $x_n$, $\eps_n$, we get that for $x=\eps_n^{\frac{3}{2}}\lam_n^{-\frac{1}{2}}|x_n-p_i|^{-\frac{\alpha}{2}} y+x_n$
\begin{eqnarray*}
\lam_n^{-\frac{1}{3}} \big( d(x)^{-\frac{\alpha}{3}} |x-x_n|^{-\frac{2}{3}}\left( 1-u_n(x) \right)\big)
=|y|^{-\frac{2}{3}} U_n(y) \big|\eps_n ^{\frac{3}{2}}\lam_n^{-\frac{1}{2}}|x_n-p_i|^{-\frac{\alpha+2}{2}} y+\frac{x_n-p_i}{|x_n-p_i|} \big|^{-\frac{\alpha}{3}}
\geq C_R>0 \end{eqnarray*}
uniformly in $n$ and $y \in B_R(0)$, for any $R>0$. Then,
$$ \frac{\eps_n^{\frac{3}{2}} \lam_n^{-\frac{1}{2}} |x_n-p_i|^{-\frac{\alpha}{2}}}{|x_n-y_n|} \to 0 \quad \hbox{as }n \to +\infty,$$
and hence, by (\ref{disgiunte}) $\phi_n$, $\psi_n$ have disjoint compact supports leading to $\mu_{2,n}<0$ in contradiction with (\ref{morse}).\\
2) if $\eps_n^{-3}\lam_n |x_n-p_i|^{\alpha+2} \leq C$ as $n \to +\infty$, given $U_n$ be as in (\ref{un3caso}) associated to $x_n$, $\eps_n$, we get that for $x=\eps_n^{\frac{3}{2+\alpha}}\lam_n^{-\frac{1}{2+\alpha}}y+x_n$
\begin{eqnarray*}
\lam_n^{-\frac{1}{3}} \big( d(x)^{-\frac{\alpha}{3}} |x-x_n|^{-\frac{2}{3}}\left( 1-u_n(x) \right)\big)
&=&|y|^{-\frac{2}{3}} U_n(y) \big|y+ \eps_n ^{-\frac{3}{2+\alpha}}\lam_n^{\frac{1}{2+\alpha}} (x_n-p_i) \big|^{-\frac{\alpha}{3}}\\
&\geq& D_R |y|^{-\frac{2}{3}} U_n(y) \geq C_R>0 \end{eqnarray*}
uniformly in $n$ and $y \in B_R(0)$, for any $R>0$. Then,
$$ \frac{\eps_n^{\frac{3}{2+\alpha}} \lam_n^{-\frac{1}{2+\alpha}}}{|x_n-y_n|} \to 0 \quad \hbox{as }n \to +\infty,$$
and hence, by (\ref{disgiunte}) $\phi_n$, $\psi_n$ have disjoint compact supports leading to a contradiction.

\medskip \noindent \medskip \noindent {\bf $3^{\hbox{rd}}$ Case} Assume that $y_n \to p_i$ in a ``fast'' way:
$$\mu_n^{-3}\lam_n |y_n-p_i|^{\alpha+2} \leq C.$$
Since $d(y_n)=|y_n-p_i|$, by (\ref{conseq}) we get that
\begin{equation} \label{conseqter} \frac{|y_n-p_i|}{|x_n-y_n|}=
\frac{\mu_n^{\frac{3}{2}} \lam_n^{-\frac{1}{2}}}{|x_n-y_n| \: |y_n-p_i|^{\frac{\alpha}{2}} }
\left( \mu_n^{-3}\lam_n |y_n-p_i|^{\alpha+2} \right)^{\frac{1}{2}} \to 0 \quad \hbox{as }n\to +\infty,\end{equation}
and then, for $n$ large:
\begin{equation} \label{conseqbis} \frac{|x_n-p_i|}{|y_n-p_i|}\geq \frac{|x_n-y_n|}{|y_n-p_i|}-1\geq 1\:\:,\:\:\:
\frac{|x_n-p_i|}{|x_n-y_n|}\geq 1- \frac{|y_n-p_i|}{|x_n-y_n|}\geq \frac{1}{2}.\end{equation}
Since $\eps_n \leq \mu_n$, by (\ref{conseq}) and (\ref{conseqbis}) we get that
\begin{eqnarray} \label{vincoio} \eps_n^{-3}\lam_n |x_n-p_i|^{\alpha+2} &\geq&  \left( \frac{\mu_n^{\frac{3}{2}} \lam_n^{-\frac{1}{2}} }{|x_n-y_n| \: |y_n-p_i|^{\frac{\alpha}{2}} } \right)^{-2} \left(\frac{|x_n-p_i|}{|x_n-y_n|^{\frac{2}{\alpha+2}} \: |y_n-p_i|^{\frac{\alpha}{\alpha+2}} } \right)^{\alpha+2}\\
&\geq& C \left( \frac{\mu_n^{\frac{3}{2}} \lam_n^{-\frac{1}{2}} }{|x_n-y_n| \: d(y_n)^{\frac{\alpha}{2}} } \right)^{-2}\to +\infty \quad \hbox{as }n \to +\infty.\nonumber
\end{eqnarray}
The meaning of (\ref{vincoio}) is the following: once $y_n$ provides a fast blowing up sequence at $p_i$, then no other fast blow up at $p_i$ can occurr as (\ref{vincoio}) states for $x_n$.

\medskip \noindent Let $\beta_n=\mu_n^{\frac{3}{2+\alpha}} \lam_n^{-\frac{1}{2+\alpha}}$. By (\ref{conseq}) and (\ref{conseqter}) we get that
\begin{eqnarray} \label{valesolo}
\frac{\beta_n}{|x_n-y_n|}=\mu_n^{\frac{3}{2+\alpha}} \lam_n^{-\frac{1}{2+\alpha}}|x_n-y_n|^{-1}=
\left( \frac{\mu_n^{\frac{3}{2}} \lam_n^{-\frac{1}{2}}}{|x_n-y_n| \: d(y_n)^{\frac{\alpha}{2}} }\right)^{\frac{2}{2+\alpha}} \left( \frac{|y_n-p_i|}{|x_n-y_n|} \right)^{\frac{\alpha}{2+\alpha}}
\to 0 \quad \hbox{as }n \to +\infty.
\end{eqnarray}
However, since $u_n$ blows up fast at $p_i$ along $y_n$, we have that $\beta_n^{-1} d(y_n)\leq C$ and then, (\ref{assRn}) does not hold. Letting as before
$$\hat U_n(y)=\frac{1-u_n(\beta_n y+y_n)}{\mu_n}\:,\:\:  y \in \hat \Omega_n=\frac{\Omega-y_n}{\beta_n},$$
we need to refine the analysis before in order to get some estimate for $\hat U_n$ even when only (\ref{valesolo}) does hold. Formula (\ref{3:18}) gives that:
\begin{eqnarray}  \hat U_n(y) &\geq&  |y_n-p_i|^{-\frac{\alpha}{3}} |y_n-x_n|^{-\frac{2}{3}} |\beta_n y+y_n-p_i|^{\frac{\alpha}{3}} |\beta_n y+y_n-x_n|^{\frac{2}{3}} \label{allora} \\
&=&
\left( \mu_n^{-\frac{3}{2+\alpha}} \lam_n^{\frac{1}{2+\alpha}} |y_n-p_i| \right)^{-\frac{\alpha}{3}}
|y+\mu_n^{-\frac{3}{2+\alpha}} \lam_n^{\frac{1}{2+\alpha}}(y_n-p_i)|^{\frac{\alpha}{3}}
|\frac{\beta_n}{|x_n-y_n|} y+\frac{y_n-x_n}{|x_n-y_n|}|^{\frac{2}{3}}\nonumber \\
&\geq & C \left( \mu_n^{-\frac{3}{2+\alpha}} \lam_n^{\frac{1}{2+\alpha}} |y_n-p_i| \right)^{-\frac{\alpha}{3}}
|y+\mu_n^{-\frac{3}{2+\alpha}} \lam_n^{\frac{1}{2+\alpha}}(y_n-p_i)|^{\frac{\alpha}{3}} \nonumber \end{eqnarray}
for $|y|\leq R_n=(\frac{|x_n-y_n|}{\beta_n})^{\frac{1}{2}}$, and $R_n \to +\infty$ as $n \to +\infty$ by (\ref{valesolo}).
Since (\ref{allora}) implies that (\ref{lontanozerobis}) holds for $\mu_n$, $y_n$, $\hat U_n$, Proposition \ref{prop3caso} provides some $\psi_n \in C_0^\infty(\Omega)$ such that (\ref{testfunction}) holds and $\hbox{Supp }\psi_n \subset B_{M \mu_n^{\frac{3}{2+\alpha}} \lam_n^{-\frac{1}{2+\alpha}}}(y_n)$, $M>0$.

\medskip \noindent Since $y_n$ cannot lie in any ball centered at $x_n$ and radius of order of the scale parameter ($\eps_n^{\frac{3}{2}}\lam_n^{-\frac{1}{2}}$ or $\eps_n^{\frac{3}{2}} \lam_n^{-\frac{1}{2}}|x_n-p_i|^{-\frac{\alpha}{2}}$), by (\ref{valesolo}) we get that $\phi_n$ and $\psi_n$ have disjoint compact supports leading to $\mu_{2,n}<0$. A contradiction to (\ref{morse}). The proof of the Proposition is now complete.\qed

\subsection{Compactness issues}
We are now in position to give the proof of Theorem \ref{thm2}. Assume $2\leq N \leq 7$. Let $f \in C(\bar \Omega)$ be as in (\ref{assf1}). Let $(\lambda_n)_n$ be a sequence such that $\lambda_n \to  \lam \in[0,\lam^*]$ and let $u_n$ be an associated solution such that (\ref{morse}) holds, namely
$$\mu_{2,n}:=\mu_{2,\lambda_n}(u_n)\geq 0.$$
The essential ingredient will be the estimate of Proposition \ref{propest} combined with the uniqueness result of Proposition \ref{uniqueness}.

\medskip \noindent {\bf Proof (of Theorem \ref{thm2}):} Let $x_n$ be the maximum point of $u_n$ in $\Omega$ and, up to a subsequence, assume by contradiction that $u_n(x_n)=\max\limits_{\Omega}u_n(x)\to 1$ as $n\to \infty $. Proposition \ref{propest} gives that:
$$u_n(x) \leq 1- C \lam _n^{\frac{1}{3}} d(x)^{\frac{\alpha}{3}} |x-x_n|^{\frac{2}{3}}$$
for any $x \in \Omega$ and $n \ge N_0$, for some $C>0$ and $N_0 \in \N$ large. Here, $d(x)=\min \{|x-p_i|:\:i=1,\dots,k \}$ stands for the distance function from the zero set of $f(x)$. Thus, we have that:
\begin{equation} \label{cheperme} 0\le \frac{\lam_n f(x)}{(1-u_n)^2}\le C \frac{f(x)}{d(x)^{\frac{2\alpha}{3}}} \frac{\lam_n^{\frac{1}{3}}}{|x-x_n|^{\frac{4}{3}}} \end{equation}
for any $x\in \Omega$ and $n\geq N_0$. Since by (\ref{assf1})
$$\big| \frac{f(x)}{d(x)^{\frac{2\alpha}{3}}}\big| \leq |x-p_i|^{\frac{\alpha}{3}} \|f_i\|_\infty \leq C$$
for $x$ close to $p_i$, $f_i$ as in (\ref{deffi}), we get that $\frac{f(x)}{d(x)^{\frac{2\alpha}{3}}}$ is a bounded function on $\Omega$ and then, by (\ref{cheperme}) $\lam_n f(x)/(1-u_n)^2$ is uniformly bounded in
$L^s(\Omega)$, for any $1<s<\frac{3N}{4}$. Standard elliptic
regularity theory now implies that $u_n$ is uniformly bounded in
$W^{2,s}(\Omega)$. By Sobolev's imbedding theorem, $u_n$ is uniformly bounded in
$C^{0,\beta}(\bar \Omega )$ for any $0<\beta <2/3$. Up to a
subsequence, we get that $u_n\to u_0$ weakly in $H_0^1(\Omega)$ and strongly in $C^{0,\beta }(\bar
\Omega )$, $0<\beta <2/3$, where $u_0$ is an H\"olderian function solving weakly in $H_0^1(\Omega)$ the equation:
\begin{equation} \left\{
\begin{array}{ll} -\Delta u_0 = \displaystyle\frac{\lambda
f(x)}{(1-u_0)^2}    & \hbox{in }\Omega\,, \\
0\leq u_0 \le 1  & \hbox{in }\Omega \,,\\
u_0=0 & \hbox{on }\partial \Omega. \end{array}\right. \label{3:23}\end{equation}
Moreover, by uniform convergence
$$\max_{\Omega}u_0=\lim_{n \to +\infty} \max_{\Omega} u_n=1$$
and, in particular $u_0>0$ in $\Omega$. Clearly, $\lam>0$ since any weak harmonic function in $H_0^1(\Omega)$ is identically zero. To reach a contradiction, we shall first show that $\mu_{1,\lam}(u_0) \geq 0$ and then deduce from the uniqueness, stated in Proposition \ref{uniqueness}, of the semi-stable solution $u_\lam$ that $u_0=u_\lam$. But $\displaystyle\max_{\Omega} u_\lam<1$ for any $\lam\in [0,\lam^*]$, contradicting $\displaystyle \max_{\Omega}u_0=1$. Hence, the claimed compactness must hold.

\medskip \noindent In addition to (\ref{morse}), assume now that  $\mu_{1,n}<0$, then $\lam>0$. Indeed, if $\lam_n \to 0$, then by compactness and standard regularity theory, we get that $u_n \to u_0$ in $C^2(\bar \Omega)$, where $u_0$ is an harmonic function so that $u_0=0$ on $\partial \Omega$. Then, $u_0=0$ and $u_n \to 0$ in $C^2(\bar \Omega)$. But the only branch of solutions for $(S)_\lam$ bifurcating from $0$ for $\lam$ small is the branch of minimal solutions $u_\lam$ and then, $u_n =u_{\lam_n}$ for $n$ large contradicting $\mu_{1,n}<0$.

\medskip \noindent  In order to complete the proof, we need only to show that
\begin{equation}
\mu_{1,\lam}(u_0)=\inf\left\{ \int_\Omega \left(|\nabla \phi |^2 -\displaystyle\frac{2\lam
f(x)}{(1-u_0)^3}\phi^2 \right) ;\,  \phi \in C_0^\infty(\Omega) \, {\rm and }\, \int_\Omega \phi^2=1 \right\}\geq 0.\label{3:24}\end{equation}

\medskip \noindent Indeed, first by Propositions \ref{prop1caso}-\ref{prop3caso} we get the existence of a function $\phi_n \in C_0^\infty(\Omega)$ so that
\begin{eqnarray}
\int_\Omega \left(|\nabla \phi _n|^2- \displaystyle\frac{2\lam_nf(x)}{(1-u_n)^3}\phi_n^2\right)<0. \label{3:26palle}
\end{eqnarray}
Moreover, $\hbox{Supp }\phi_n \subset B_{r_n}(x_n)$ and $r_n \to 0$ as $n\to +\infty$. Up to a subsequence, assume that $x_n \to p \in \bar \Omega$ as $n \to +\infty$.

\medskip \noindent By contradiction, if (\ref{3:24}) were false, then
there exists $\phi_0 \in C^{\infty }_0(\Omega )$ such that
\begin{equation}
\int_\Omega \left(|\nabla \phi _0|^2 -\displaystyle\frac{2\lam
f(x)}{(1-u_0)^3}\phi_0^2 \right)  < 0. \label{3:25} \end{equation}
We will replace $\phi_0$ with a truncated function $\phi_\delta$ with $\delta>0$ small enough, and so that (\ref{3:25}) is still true while $\phi_\delta=0$ in $B_{\delta^2}(p) \cap \Omega$. In this way, $\phi_n$ and $\phi_\delta$ would have disjoint compact supports in contradiction to $\mu_{2,n}\geq 0$.

\medskip \noindent Let $\delta>0$ and set $\phi_\delta=\chi_\delta \phi_0 $, where $\chi_\delta$ is a cut-off function defined as:
$$\chi_\delta(x)=
\begin{cases}0 \quad & |x-p|\le \delta^2\,,\\[1mm]
2-\displaystyle\frac{\log |x-p|}{\log \delta} \,\quad & \delta^2 \le  |x-p|\le \delta\,,\\[1mm]
1\,\quad &  |x-p|\ge \delta\,.
\end{cases}$$
By Lebesgue's theorem, we have:
\begin{equation}
\int_\Omega \frac{2\lam f(x)}{(1-u_0)^3}\phi_\delta^2\to \int_\Omega \frac{2\lam
f(x)}{(1-u_0)^3}\phi _0^2\,,\quad \hbox{as } \delta \to 0\,.
\label{3:25''}\end{equation}
For the gradient term, we have the expansion:
$$\int_\Omega |\nabla \phi_\delta|^2=\int_\Omega \phi_0^2 |\nabla \chi_\delta|^2+\int_\Omega \chi_\delta^2 |\nabla \phi_0|^2+2\int_\Omega \chi_\delta \phi _0 \nabla \chi_\delta  \nabla \phi _0\,.$$
The following estimates hold:
$$0\leq \int_\Omega \phi_0^2 |\nabla \chi_\delta|^2 \leq \|\phi_0\|_\infty^2 \displaystyle\int_{\delta^2 \le |x-p|\le \delta}
\displaystyle\frac{1}{|x-p|^2\log ^2 \delta} \le \displaystyle \frac{C}{\log \frac{1}{\delta}}$$
and
$$\big | 2\displaystyle\int_\Omega \chi _\delta \phi _0 \nabla \chi_\delta \nabla \phi _0 \big|
\leq  \displaystyle\frac{2\|\phi_0\|_\infty  \|\nabla \phi_0\|_{\infty}}{\log \frac{1}{\delta}} \displaystyle\int_{B_1(0)}\displaystyle\frac{1}{|x|},$$
and provide:
\begin{equation}\int_\Omega |\nabla \phi_\delta|^2 \to \int_\Omega |\nabla \phi_0|^2 \quad  \hbox{as } \delta \to 0.\label{3:26}\end{equation}
Combining (\ref{3:25})-(\ref{3:26}), we get that:
$$\int_\Omega \left(|\nabla \phi_\delta|^2-\frac{2\lam f(x)}{(1-u_0)^3}\phi_\delta^2 \right) <0$$
for $\delta>0$ sufficiently small. This completes the proof of (\ref{3:24}) and Theorem \ref{thm2} is completely established.

\section{The one dimensional problem}
Let $I=(a,b)$ be a bounded interval in $\R$. Assume $f \in C^1(\bar I)$ so that $f \geq C>0$ in $I$. In Theorem \ref{thm5} we study solutions $u_n$ of the following problem:
\begin{equation}
\left\{ \begin{array}{ll} -  \ddot u_n = \displaystyle\frac{\lambda _n
f(x)}{(1-u_n)^2}  & \hbox{in }I\,,\\
0<u_n<1  & \hbox{in } I \,,\\
u_n(a)=u_n(b)=0. & \end{array} \right. \label{3:33}
\end{equation}

\medskip \noindent {\bf Proof (of Theorem \ref{thm5}):} Assume that $u_n$ satisfy (\ref{morsebound}) and $\lam_n \to \lambda \in (0,\lam^*]$. Let $x_n \in I$ be a maximum point: $u_n(x_n)=\displaystyle \max_{I}u_n$. If $(u_n)_n$ is not compact, then  up to a subsequence, we may assume that $u_n(x_n) \to 1$ with $x_n \to x_0 \in \bar I$ as $n \to +\infty$. Away from $x_0$, $u_n$ is uniformly far away from $1$. Otherwise, by the maximum principle we would have $u_n\to 1$ on an interval of positive measure, and then $\mu_{k,\lam_n}(u_n)<0$, for any $k$ and $n$ large. A contradiction.

\medskip \noindent Assume, for example, that $a\leq x_0 <b$. By elliptic regularity theory, $\dot u_n$ is uniformly bounded far away from $x_0$. Let $\eps>0$. We multiply (\ref{3:33}) by $\dot u_n$ and integrate on $(x_n,x_0+\eps)$:
$$\dot u_n^2(x_n)- \dot u_n^2(x_0+\eps) = \int_{x_n}^{x_0+\eps} \frac{2 \lam_nf(s) \dot u_n(s)}
{(1-u_n(s))^2}ds= \frac{2 \lam_nf(x_0+\eps)}{1-u_n(x_0+\eps)}- \frac{2 \lam_nf(x_n)}{1-u_n(x_n)}
- \int_{x_n}^{x_0+\eps} \frac{2 \lam_n\dot f(s)}{1-u_n(s)}ds.$$
Then, for $n$ large:
\begin{eqnarray*} \dot u_n^2(x_n)+ \frac{C \lam }{1-u_n(x_n)} &\leq& \dot u_n^2(x_0+\eps)+ 2 \lam_n \frac{f(x_0+\eps)}{1-u_n(x_0+\eps)}
-2 \lam_n \int_{x_n}^{x_0+\eps} \frac{\dot f(s)}{1-u_n(s)}ds\\
&\leq& C_\eps+ 4 \lam \parallel \dot f \parallel_\infty \frac{x_0+\eps-x_n}{1-u_n(x_n)}
\end{eqnarray*}
since $u_n(x_n)$ is the maximum value of $u_n$ in $I$. Choosing $\eps>0$ sufficiently small, we get that
for any $n$ large: $\frac{1}{1-u_n(x_n)}\leq C_\eps$, contradicting $u_n(x_n) \to 1$ as $n \to +\infty$. \qed

\section{The second bifurcation point and the branch of unstable solutions}
We now establish Theorem \ref{thm4}. First, let us recall the definition of $\lambda_2^*$:
$$\lambda_2^*=\inf \{ \beta>0:\exists \hbox{ a curve}  V_\lambda \in C([\beta ,\lambda^*]; C^2(\Omega)) \hbox{ of solutions to } (S)_\lambda \hbox{ s.t. }\mu_{2,\lambda}(V_\lambda) \geq 0  , \:V_\lambda \equiv U_\lambda \: \forall \lambda \in (\lambda^*-\delta,\lambda^*) \}.$$
As for as Theorem \ref{thm4} is concerned, for any $\lam \in (\lam_2^*,\lam^*)$ by definition there exists a solution $V_\lam$ and it is such that:
\begin{equation} \label{aloms}\mu_{1,\lam}:=\mu_{1,\lam}(V_\lam)<0 \quad \forall\: \lam \in (\lambda_2^*,\lam^*).\end{equation}
In particular, $V_\lam \not= u_\lam$ provides a second solution different from the minimal one.

\medskip \noindent Clearly (\ref{aloms}) is true because first $\mu_{1,\lam}<0$ for $\lambda$ close to $\lambda^*$. Moreover, if  $\mu_{1,\lam}=0$ for some $\lam \in (\lam_2^*,\lam^*)$, then by Proposition \ref{uniqueness} $V_\lam=u_\lam$ contradicting the fact that $\mu_{1,\lam}(u_\lam)>0$ for any $0<\lam<\lam^*$.

\medskip \noindent Since by definition $\mu_{2,\lam}(V_\lam)\geq 0$ for any $\lam \in (\lam_2^*,\lam^*)$, we can take a sequence $\lam_n \downarrow \lam_2^*$ and apply Theorem \ref{thm2} to get that $\lam_2^*=\displaystyle \lim_{n \to +\infty} \lam_n>0$, $\displaystyle \sup_{n \in \N}\parallel V_{\lam_n} \parallel_\infty <1$. By elliptic regularity theory, up to a subsequence $V_{\lam_n} \to V^*$ in $C^2(\bar \Omega)$, where $V^*$ is a solution for $(S)_{\lam_2^*}$. As before, $\mu_{1,\lam_2^*}(V^*)<0$ and by continuity $\mu_{2,\lam_2^*}(V^*)\geq 0$.

\medskip \noindent If $\mu_{2,\lam_2^*}(V^*)>0$, let us fix some $\eps>0$ small so that $0\leq V^* \leq 1-2 \eps$ and consider the truncated nonlinearity $g_\eps(u)$ as in (\ref{geps}). Clearly, $V^*$ is a solution of (\ref{2:1}) at $\lam=\lam_2^*$ so that $-\Delta-\lam_2^* f(x) g_\eps'(V^*)$ has no zero eigenvalues. Namely, $V^*$ solves $N(\lam_2^*, V^*)=0$, where $N$ is a map from $\R \times C^{2,\alpha}(\bar \Omega)$ into $C^{2,\alpha}(\bar \Omega)$, $\alpha \in (0,1)$, defined as:
$$N: (\lam,V) \:\: \longrightarrow \:\: V+\Delta^{-1}\left( \lam f(x)g_\eps(V)\right).$$
Moreover,
$$\partial_V N(\lam_2^*,V^*)=\hbox{Id }+\Delta^{-1}\left( \displaystyle \frac{2 \lam_2^* f(x)}{(1-V^*)^3}\right)$$
is an invertible map since $-\Delta-\lam_2^* f(x) g_\eps'(V^*)$ has no zero eigenvalues. The Implicit Function Theorem gives the existence of a curve $W_\lam$, $\lam \in (\lam_2^*-\delta,\lam_2^*+\delta)$, of solution for (\ref{2:1}) so that $\lim_{\lam \to \lam_2^*}W_\lam=V^*$ in $C^{2,\alpha}(\bar \Omega)$. Up to take $\delta$ smaller, this convergence implies that $\mu_{2,\lam}(W_\lam)>0$ and $W_\lam\leq 1-\eps$ for any $\lam \in (\lam_2^*-\delta,\lam_2^*+\delta)$. Hence, $W_\lam$ is a solution of $(S)_\lam$ so that $\mu_{2,\lam}(W_\lam)>0$ contradicting the definition of $\lam_2^*$. Hence, $\mu_{2,\lam_2^*}(V^*)=0$. A similar argument works for the radial problem $(S)_\lam$ on the unit ball and $f(x)$ as in (\ref{assf0}), provided either $2\leq N \leq7$ or $N\geq 8$, $\alpha>\alpha_N$. The proof of Theorem \ref{thm4} is complete.

\section{Appendix}
We shall prove here the following Theorem already announced in the Introduction.
\begin{thm} \label{thm3.0} Assume either $1\le N\le 7$ or $N\geq 8$, $\alpha>\alpha_N$. Let $U$ be a solution of
\begin{equation} \left\{ \begin{array}{ll} \Delta U = \displaystyle\frac{|y|^\alpha}{U^2}
& \hbox{in }\R ^N,\\
U(y) \geq C>0  & \hbox{in } \R^N. \end{array} \right.  \label{3:1agg.bis}
\end{equation}
Then,
\begin{equation}
\mu _1(U)=\inf \left\{ \int_{\R^N} \big( |\nabla \phi
|^2- \frac{2 |y|^\alpha}{U^3} \phi ^2\big); \, \phi \in C^{\infty }_0(\R^N)\, {\rm and}\,  \int_{\R^N} \phi^2=1 \right\}<0\,. \label{3:2agg.bis}
\end{equation}
Moreover, if $N\geq 8$ and $0\leq \alpha\leq \alpha_N$, then there exists at least a  solution $U$ of (\ref{3:1agg.bis}) such that $\mu_1(U) \geq 0$.
\end{thm}
\noindent{\bf Proof:} By contradiction,  assume that
$$\mu _1(U)=\inf \left\{ \int_{\R^N} \big( |\nabla \phi
|^2- \frac{2 |y|^\alpha}{U^3} \phi ^2\big)\, ; \, \phi \in C^{\infty }_0(\R^N)\, {\rm and}\,  \int_{\R^N} \phi^2\, dx=1 \right\} \geq 0.$$
By the density of $C_0^\infty(\R^N)$ in $D^{1,2}(\R^N)$, we have that
\begin{equation}
\int |\nabla \phi |^2\ge 2 \int \frac{|y|^\alpha}{U^3}\phi ^2\,, \quad \forall \:
\phi \in D^{1,2}(\R^N)\,. \label{3:3}
\end{equation}
In particular, the test function $\phi=\frac{1}{(1+|y|^2)^{\frac{N-2}{4}+\frac{\delta
}{2}}}\in D^{1,2}(\R^N)$ applied in  (\ref{3:3}) gives that
\begin{equation}
\int \displaystyle\frac{|y|^\alpha }{(1+|y|^2)^{\frac{N-2}{2}+\delta
}U^{3}} \le C \int \frac{1}{(1+|y|^2)^{\frac{N}{2}+\delta }}<+\infty
\,,\label{3:7}\end{equation}
for any $\delta>0$.

\vspace {0.5cm} \noindent {\bf Step 1.}\ \ We want to show that (\ref{3:3}) allows us to perform the following Moser-type iteration scheme: for any $0<q<4+2\sqrt 6$ and $\beta$ there holds
\begin{equation} \displaystyle \int\displaystyle\frac{1}{(1+|y|^2)^{\beta
-1-\frac{\alpha}{2} }U^{q+3}}\le C_q \left(1+ \displaystyle\int\displaystyle\frac{1}{(1+|y|^2)^{\beta
}U^q} \right) \,\label{3:6} \end{equation}
(provided the second integral is finite).

\medskip \noindent Indeed, let $R>0$ and consider a smooth radial cut-off
function $\eta$ so that: $0\leq \eta \leq 1$, $\eta =1$ in $B_R(0)$, $\eta=0$ in
$\R^N \setminus B_{2R}(0)$. Multiplying (\ref{3:1agg.bis}) by $\displaystyle \frac{\eta ^2}{(1+|y|^2)^{\beta
-1}U^{q+1}}$, $q>0$, and integrating by parts we get:
\begin{eqnarray*}
\displaystyle\int \displaystyle\frac{|y|^\alpha \eta ^2 }{(1+|y|^2)^{\beta
-1}U^{q+3}}&=& \displaystyle\frac{4(q+1)}{q^2}\displaystyle\int\Big |\nabla
\Big(\displaystyle\frac{\eta }{(1+|y|^2)^{\frac{\beta
-1}{2}}U^{\frac{q}{2}}}\Big)\Big
|^2-\displaystyle\frac{4(q+1)}{q^2}\displaystyle\int\displaystyle\frac{1}{U^q}\Big
|\nabla \Big(\displaystyle\frac{\eta }{(1+|y|^2)^{\frac{\beta
-1}{2}}}\Big)\Big
|^2\\
&&-\displaystyle\frac{q+2}{q^2} \displaystyle\int\displaystyle \nabla(\frac{1}{U^q}) \nabla \Big(\displaystyle\frac{\eta ^2}{(1+|y|^2)^{\beta
-1}}\Big) \\
&=&\displaystyle\frac{4(q+1)}{q^2}\displaystyle\int\Big |\nabla
\Big(\displaystyle\frac{\eta }{(1+|y|^2)^{\frac{\beta
-1}{2}}U^{\frac{q}{2}}}\Big)\Big
|^2-\displaystyle\frac{2}{q}\displaystyle\int\displaystyle\frac{1}{U^q}\Big
|\nabla \Big(\displaystyle\frac{\eta }{(1+|y|^2)^{\frac{\beta
-1}{2}}}\Big)\Big
|^2\\
&&+\displaystyle\frac{2(q+2)}{q^2}\displaystyle\int\displaystyle\frac{1}{U^q}\displaystyle\frac{\eta
}{(1+|y|^2)^{\frac{\beta -1}{2}}}\Delta \Big(\displaystyle\frac{\eta
}{(1+|y|^2)^{\frac{\beta -1}{2}}}\Big)\,,
\end{eqnarray*}
by means of the relation: $\Delta (\psi)^2=2|\nabla \psi|^2+2 \psi \Delta \psi$.

\medskip \noindent Then, by (\ref{3:3}) we deduce that
$$(8q+8-q^2)\displaystyle\int\frac{|y|^\alpha \eta ^2 }{(1+|y|^2)^{\beta-1}U^{q+3}}\le C_q'\int \frac{1}{U^q}\Big( \big |\nabla \big(\displaystyle\frac{\eta }{(1+|y|^2)^{\frac{\beta-1}{2}}}\big)\big |^2+\displaystyle\frac{\eta
}{(1+|y|^2)^{\frac{\beta -1}{2}}} \big|\Delta \big(\displaystyle\frac{\eta
}{(1+|y|^2)^{\frac{\beta -1}{2}}}\big)\big| \Big).$$
Assuming that $|\nabla \eta|\leq \frac{C}{R}$ and $|\Delta \eta|\leq \frac{C}{R^2}$, it is straightforward to see that:
$$\big |\nabla
\big(\displaystyle\frac{\eta }{(1+|y|^2)^{\frac{\beta-1}{2}}}\big)\big |^2+\displaystyle\frac{\eta
}{(1+|y|^2)^{\frac{\beta-1}{2}}} \big|\Delta \big(\displaystyle\frac{\eta
}{(1+|y|^2)^{\frac{\beta -1}{2}}}\big)\big| \leq C
\Big(\frac{1}{(1+|y|^2)^{\beta}}+\displaystyle\frac{1}{R^2(1+|y|^2)^{\beta -1}}\chi_{B_{2R}(0) \setminus B_R(0) }\Big]$$
for some constant $C$ independent on $R>0$. Then,
$$(8q+8-q^2)\displaystyle\int\frac{|y|^\alpha \eta ^2 }{(1+|y|^2)^{\beta
-1}U^{q+3}}\le C_q'' \int \frac{1}{(1+|y|^2)^{\beta} U^q}.$$

\medskip \noindent Let $q_{+}=4+2\sqrt 6$. For any $0<q<q_+$, we have $8q+8-q^2>0$ and therefore:
$$\displaystyle\int\displaystyle\frac{|y|^\alpha \eta ^2 }{(1+|y|^2)^{\beta
-1}U^{q+3}}\le C_q \displaystyle\int\displaystyle\frac{1}{(1+|y|^2)^{\beta}U^q}\,,$$
where $C_q$ does not depend on $R>0$. Taking the limit as $R \to +\infty$, we get that:
$$\displaystyle\int\displaystyle\frac{|y|^\alpha}{(1+|y|^2)^{\beta
-1}U^{q+3}}\le C_q \displaystyle\int\displaystyle\frac{1}{(1+|y|^2)^{\beta}U^q}\,$$
and then, the validity of (\ref{3:6}) easily follows.

\vspace {0.5cm} \noindent {\bf Step 2.}\ \ Let now $1\leq N \leq 7$ or $N\geq 8$, $\alpha>\alpha_N$. We want to show that
\begin{equation}
\int\displaystyle\frac{1}{(1+|y|^2)U^q}<+\infty  \label{3:13}
\end{equation}
for some $0<q<q_+=4 +2\sqrt{6}.$

\medskip \noindent Indeed, set $\beta_0=\frac{N-2-\alpha}{2}+\delta$, $\delta>0$, and $q_0=3$. By (\ref{3:7}) we get that
$$\int \displaystyle\frac{1}{(1+|y|^2)^{\beta_0}U^{q_0}} <+\infty.$$
Let $\beta_i=\beta_0-i(1+\frac{\alpha}{2})$ and $q_i=q_0+3i$, $i \in \N$. Since $q_0<q_1<q_+=4+2\sqrt{6}<q_2$, we can iterate (\ref{3:6}) exactly two times to get that:
\begin{equation} \label{itter} \int \displaystyle\frac{1}{(1+|y|^2)^{\beta_2}U^{q_2}} <+\infty \end{equation}
where $\beta_2=\frac{N-6-3\alpha}{2}+\delta$, $q_2=9$.

\medskip \noindent \medskip \noindent Let $0<q<q_+=4+2\sqrt{6}< 9$. By (\ref{itter}) and H\"older inequality we get that:
\begin{eqnarray*}
\displaystyle\int\displaystyle\frac{1}{(1+|y|^2)U^q}&=&
\displaystyle\int\displaystyle\frac{(1+|y|^2)^{\frac{q}{9}(\frac{6-N}{2}-\delta+\frac{3}{2}\alpha
)}}{U^q}\cdot\displaystyle\frac{1}{(1+|y|^2)^{\frac{q}{9}(\frac{6-N}{2}-\delta+\frac{3}{2}\alpha)+1}}\\
&\le & \Big(
\displaystyle\int\displaystyle\frac{1}{(1+|y|^2)^{\beta_2} U^{q_2}}\Big)^{\frac{q}{9}}
\Big(\displaystyle\int\displaystyle\frac{1}{(1+|y|^2)^{\frac{q}{9-q}(\frac{6-N}{2}-\delta+\frac{3}{2}\alpha
)+\frac{9}{9-q}}}\Big )^{\frac{9-q}{9}}<+\infty
\end{eqnarray*}
provided $-\frac{2q}{9-q}\beta_2+\frac{18}{9-q}>N$ or equivalently
\begin{equation} \label{choice} q>\frac{9N-18}{6-2\delta+3\alpha}. \end{equation}
To have (\ref{choice}) for some $\delta>0$ small and $q_<q_+$ at the same time, we need to require $\frac{3N-6}{2+\alpha}<q_+$ or equivalently
$$1\leq N \leq 7 \quad \hbox{ or } \quad N\geq 8\:,\:\:\alpha >\alpha_N=\frac{3N-14-4\sqrt{6}}{4+2\sqrt{6}}.$$
Our assumptions then provide the existence of some $0<q<q_+=4+2\sqrt{6}$ such that (\ref{3:13}) holds.

\vspace{0.5 cm} \noindent {\bf Step 3.}\ \ We are ready to obtain a
contradiction. Let $0<q<4+2\sqrt{6}$ be such that (\ref{3:13}) holds. Let $\eta$ be the cut-off function of Step 1. Using equation (\ref{3:1agg.bis}) we compute:
\begin{eqnarray*}
\displaystyle\int \big |\nabla
\big(\displaystyle\frac{\eta}{U^{\frac{q}{2}}}\big) \big
|^2-\displaystyle\int\displaystyle\frac{2|y|^\alpha}{U^3} \left( \frac{\eta}{U^{\frac{q}{2}}} \right) &=&\displaystyle\frac{q^2}{4}\displaystyle\int\displaystyle\frac{\eta
^2|\nabla U|^2}{U^{q+2}}+\displaystyle\int\displaystyle\frac{|\nabla
\eta |^2}{U^q}+\displaystyle\frac{1}{2}\displaystyle\int \nabla (\eta^2) \nabla
\big(\displaystyle\frac{1}{U^q}\big)-\displaystyle\int\displaystyle\frac{2 |y|^\alpha \eta
^2}{U^{q+3}}\\
&=&-\displaystyle\frac{q^2}{4(q+1)}\displaystyle\int\nabla
U\cdot\nabla\big(\displaystyle\frac{\eta
^2}{U^{q+1}}\big)+\displaystyle\int\displaystyle\frac{|\nabla \eta
|^2}{U^q}\\
&&+\frac{q+2}{4(q+1)} \displaystyle\int \nabla (\eta^2) \nabla \big(\displaystyle\frac{1}{U^q}\big)
- \displaystyle\int\displaystyle\frac{2 |y|^\alpha \eta^2}{U^{q+3}}\\
&=&-\displaystyle\frac{8q+8-q^2}{4(q+1)}\displaystyle\int\displaystyle\frac{|y|^\alpha \eta
^2}{U^{q+3}}+\displaystyle\int\displaystyle\frac{|\nabla \eta
|^2}{U^q}-\displaystyle\frac{q+2}{4(q+1)}\displaystyle\int
\displaystyle\frac{\Delta \eta ^2}{U^q}\,.
\end{eqnarray*}
Since $0<q<4+2\sqrt{6}$, $8q+8-q^2>0$ and
\begin{eqnarray*}
\displaystyle\int \big |\nabla
\big(\displaystyle\frac{\eta}{U^{\frac{q}{2}}}\big) \big
|^2-\displaystyle\int\displaystyle\frac{2 |y|^\alpha}{U^3} \left(\frac{\eta}{U^{\frac{q}{2}}}\right)^2 &\le
& -\displaystyle\frac{8q+8-q^2}{4(q+1)}\displaystyle\int_{B_1(0)}\displaystyle\frac{|y|^\alpha \eta^2}{U^{q+3}}
+O\big(\displaystyle\frac{1}{R^2}\displaystyle\int_{B_{2R}(0) \setminus B_R(0) }\displaystyle\frac{1}{U^q} \big)\\
&\le &-\displaystyle\frac{8q+8-q^2}{4(q+1)}\displaystyle\int_{B_1(0)}\displaystyle\frac{|y|^\alpha \eta^2}{U^{q+3}}+O\big(
\displaystyle\int_{|y|\geq R }\displaystyle\frac{1}{(1+|y|^2)U^q}\big).
\end{eqnarray*}
Since (\ref{3:13}) implies: $\lim_{R \to +\infty} \displaystyle\int_{|y| \geq R}\displaystyle\frac{1}{(1+|y|^2)U^q}=0$, we get that for $R$ large
$$\displaystyle\int \big |\nabla
\big(\displaystyle\frac{\eta}{U^{q/2}}\big) \big
|^2-\displaystyle\int\displaystyle\frac{2 |y|^\alpha}{U^3} \left(\frac{\eta}{U^{q/2}}\right)^2 \leq -\displaystyle\frac{8q+8-q^2}{4(q+1)}\displaystyle\int_{B_1(0)}\displaystyle\frac{|y|^\alpha}{U^{q+3}}+O\big(
\displaystyle\int_{|y|\geq R }\displaystyle\frac{1}{(1+|y|^2)U^q}\big)<0.$$
A contradiction to (\ref{3:3}). Hence, (\ref{3:2agg.bis}) holds and the proof of the first part of Theorem \ref{thm3.0} is complete. \qed

\medskip \noindent To describe the counterexample, we want to compute explicitly $u^*$ and $\lam^*$ on the unit ball with $f(x)=|x|^\alpha$ and $N\geq8$, $0\leq \alpha \leq \alpha_N$. This will then provide an example of an extremal function $u^*$ which satisfies $\parallel u^* \parallel_\infty=1$ and is therefore not a classical solution. The second part of Theorem \ref{thm3.0} then follows by
considering the limit profile around zero as $\lambda \to \lambda^*$ for the minimal solution $u_\lambda$ for $(S)_\lambda$ on the unit ball with $f(x)=|x|^\alpha$.

\medskip \noindent We shall borrow ideas from \cite{BCMR,BV}, where the authors deal with the case of regular nonlinearities. However, unlike these papers where solutions are considered in a very weak sense,  we consider here a more focussed and much simpler situation. Our example is based on the following useful characterization of the extremal solution:
\begin{thm} Let $f \in C(\bar \Omega)$ be a nonnegative function. For $\lam >0$, consider $u\in H^1_0(\Omega )$ to be a weak solution of $(S)_\lam $ (in the $H_0^1(\Omega)$-sense) such  that $\parallel u \parallel_{L^\infty(\Omega)} =1$. Then
the following assertions are equivalent:
\begin{enumerate}
\item $u$ satisfies
\begin{equation}\int _\Omega |\nabla \phi |^2 \ge \int _\Omega \frac{2\lam f(x)}{(1-u)^3}\phi ^2\quad \forall \phi
\in H^1_0(\Omega )\,,\label{6:16}\end{equation}
\item $\lam =\lam ^*$ and $u=u^*$.
\end{enumerate}
\label{thmBV} \end{thm}

\medskip \noindent Here and in the sequel, $u$ will be called a $H_0^1(\Omega)$-weak solution of $(S)_\lam$ if $0\leq u\leq 1$ a.e. while $u$ solves $(S)_\lam$ in the weak sense of $H_0^1(\Omega)$. We need the following uniqueness result:
\begin{prop} \label{uniqueness} Let $f \in C(\bar \Omega)$ be a nonnegative function. Let $u_1$, $u_2$ be two $H_0^1(\Omega)$-weak solutions of $(S)_\lam$ so that $\mu_{1,\lam}(u_i)\geq 0$, $i=1,2$. Then, $u_1=u_2$ a.e. in $\Omega$. \end{prop}
\proof For any $\theta \in [0,1]$ and $\phi \in H^1_0(\Omega)$, $\phi \ge 0$, we have that:
\begin{eqnarray*}
I_{\theta ,\phi }:&=&\displaystyle\int_\Omega \nabla \big(\theta
u_1+(1-\theta )u_2 \big) \nabla \phi -\displaystyle\int_\Omega
\displaystyle\frac{\lam f(x)}{\big(1-\theta u_1-(1-\theta )u_2
\big)^2}\phi \\
&=&\lam \displaystyle\int_\Omega f(x) \Big(\displaystyle\frac{\theta
}{(1-u_1)^2}+\displaystyle\frac{1-\theta }{(1-u_2
)^2}-\displaystyle\frac{1}{\big(1-\theta u_1-(1-\theta )u_2
\big)^2}\Big) \phi \ge 0 \end{eqnarray*}
due to the convexity of $1/(1-u)$ with respect to $u$. Since $I_{0,\phi}=I_{1,\phi }=0$, the derivative of $I_{\theta,\phi }$ at $\theta=0,1$ provides:
\begin{eqnarray*}
&&\int_\Omega \nabla \big( u_1-u_2 \big)\nabla \phi -\displaystyle\int_\Omega
\displaystyle\frac{2\lam f(x)}{\big(1- u_2 \big)^3}(u_1-u_2
)\phi \ge 0\\
&& \int_\Omega \nabla \big( u_1-u_2 \big) \nabla \phi -\displaystyle\int_\Omega
\displaystyle\frac{2\lam f(x)}{\big(1- u_1 \big)^3}(u_1-u_2 )\phi
\le 0 \end{eqnarray*}
for any $\phi \in H^1_0(\Omega )$, $\phi \ge 0$.

\medskip \noindent Testing the first inequality on $\phi=(u_1-u_2)^-$ and the second one on $(u_1-u_2)^+$ we get that:
\begin{eqnarray*}
&& \int_\Omega \Big[|\nabla (u_1-u_2)^-|^2-\frac{2\lam f(x)}{(1-u_2
)^3}\big((u_1-u_2) ^-\big)^2\Big]\le 0\,\\
&& \int_\Omega \Big[|\nabla (u_1-u_2)^+|^2-\frac{2\lam f(x)}{(1-u_1
)^3}\big((u_1-u_2)^+ \big)^2\Big]\le 0\,.
\end{eqnarray*}
Since $\mu_{1,\lam}(u_1)\geq 0$, we have that:\\
(1). if $\mu_{1,\lam}(u_1)> 0$, then $u_1 \le u_2$ a.e.;\\
(2). if $\mu_{1,\lam}(u_1)= 0$, then
\begin{equation}\label{iseigen} \int_\Omega \nabla  \big( u_1-u_2 \big) \nabla \bar\phi
-\displaystyle\int_\Omega \displaystyle\frac{2\lam f(x)}{\big(1- u_1
\big)^3}(u_1-u_2)\bar \phi =0 \end{equation}
where $\bar \phi=(u_1-u_2)^+$. Since $I_{\theta ,\bar\phi }\ge 0$ for any $\theta \in [0,1]$ and $I_{1
,\bar\phi }=\partial_\theta I_{1 ,\bar\phi }=0$,
we get that:
$$\partial^2_{\theta \theta} I_{1 ,\bar\phi}=-\int _{\Omega }\frac{6 \lam f(x)}{(1-u_1)^4}\big((u_1-u_2
)^+\big)^3\ge 0\,.$$
Let $Z_0=\{x \in \Omega: \: f(x)=0 \}$. Clearly, $(u_1-u_2)^+=0$ a.e. in $\Omega \setminus Z_0$ and, by (\ref{iseigen}) we get:
$$\int_\Omega |\nabla  \big( u_1-u_2 \big)^+|^2=0.$$
Hence, $u_1 \leq u_2$ a.e. in $\Omega$. The same argument applies to prove the reversed inequality: $u_2 \leq u_1$ a.e. in $\Omega$. Therefore, $u_1=u_2$ a.e. in $\Omega$ and the proof is complete.\qed

\medskip \noindent Since $\parallel u_\lam\parallel<1$ for any $\lam \in (0,\lam^*)$,  we need --in order to prove Theorem \ref{thmBV}-- only to show that $(S)_\lam$ does not have any $H_0^1(\Omega)$-weak solution for $\lam >\lam^*$. By the definition of $\lam^*$, this is already true for classical solutions. We shall now extend this property to the class of weak solutions by means of the following result:
\begin{prop}If $w$ is a $H_0^1(\Omega)$-weak solution of $(S)_{\lam }$, then for any $\varepsilon \in (0,1)$ there exists a classic solution $w_{\varepsilon}$ of $(S)_{\lam (1-\varepsilon )}$.
\end{prop}
\proof First of all, we prove that: for any $\psi \in C^2([0,1])$  concave function so that $\psi (0)=0$, we have that
\begin{equation}
\int _{\Omega }\nabla \psi (w)  \nabla \varphi\ge \int_{\Omega } \frac{\lam f}{(1-w)^2} \dot \psi (w)\varphi  \label{supsol}
\end{equation}
for any $\varphi \in H^1_0(\Omega )$, $\varphi \ge 0$. Indeed, by concavity of $\psi$ we get:
\begin{eqnarray*} \displaystyle\int _{\Omega }\nabla \psi (w)  \nabla
\varphi &=& \displaystyle\int _{\Omega }\dot \psi(w) \nabla w \nabla \varphi =
\displaystyle\int _{\Omega } \nabla w \nabla \left( \dot \psi(w) \varphi \right)
- \displaystyle\int _{\Omega }\ddot \psi(w)\varphi  |\nabla w|^2  \\
&\ge & \displaystyle\int _{\Omega } \displaystyle\frac{\lam
f(x)}{(1-w)^2} \dot \psi (w )\varphi\, \end{eqnarray*}
for any $\varphi \in C_0^\infty(\Omega)$, $\varphi \geq 0$. By density, we get (\ref{supsol}).

\medskip \noindent Let $\eps\in (0,1)$. Define
$$\psi_\eps (w):=1-\left(\varepsilon+(1-\varepsilon )(1-w)^3\right)^{\frac{1}{3}}\,,\quad 0 \le w \leq 1\,.$$
Since $\psi_\eps \in C^2([0,1])$ is a concave function, $\psi_\eps(0)=0$ and
$$\dot \psi_\eps (w)=(1-\eps) \frac{g \left(\psi_\eps(w) \right) }{ g(w) }\:,\quad g(s):=(1-s)^{-2},$$
by (\ref{supsol}) we obtain that for any $\varphi \in H^1_0(\Omega )$, $\varphi \ge 0$:
\begin{eqnarray*}
\displaystyle\int _{\Omega }\nabla \psi_\eps (w)  \nabla \varphi\ge
\displaystyle\int _{\Omega } \displaystyle\frac{\lam  f(x)}{(1-w)^2}
\dot \psi_\eps (w) \varphi = \lam(1-\eps) \displaystyle\int _{\Omega
} f(x) g \big(\psi_\eps (w)\big) \varphi =\displaystyle\int _{\Omega
}\frac{\lam (1-\varepsilon )f(x)}{(1-\psi_\eps(w))^2}\varphi\,.
\end{eqnarray*} Hence, $\psi_\eps(w)$ is a $H_0^1(\Omega)$-weak
supersolution of $(S)_{\lam(1-\eps)}$ so that $0 \leq
\psi_\eps(w)\leq 1-\eps^{\frac{1}{3}}<1$. Since $0$ is a subsolution
for any $\lam>0$, we get the existence of a $H_0^1(\Omega)$-weak
solution $w_\eps$ of $(S)_{\lam(1-\eps)}$ so that $0\leq w_\eps \leq
1-\eps^{\frac{1}{3}}$. By standard elliptic regularity theory,
$w_\eps$ is a classical solution of $(S)_{\lam(1-\eps)}$. \qed

\bigskip \noindent We are now ready to provide the counterexample on $B=B_1(0)$. We want to show that $u^*(x)=1-|x|^{\frac{2+\alp }{3}}$ and $\lam_*=\frac{(2+\alp )(3N+\alp -4)}{9}$. It is easy to check that $u^*$ is a $H_0^1(\Omega)$-weak solution of $(S)_{\lam^*}$, provided $\alpha>1$ if $N=1$ and $\alpha \geq 0$ if $N\geq 2$. By the characterization of Theorem \ref{thmBV}, we need only to prove (\ref{6:16}). By Hardy's inequality, we have that for $N\geq 2$:
$$\int _B |\nabla \phi |^2 \geq  \frac{(N-2)^2}{4} \int_B \frac{ \phi^2}{|x|^2}$$
for any $\phi \in H^1_0(B)$, and then (\ref{6:16}) holds if $2\lam^*\leq \frac{(N-2)^2}{4}$, or equivalently, if
$$ N \ge 8 \quad \hbox{and} \quad 0 \le \alp \le  \alp_N.$$

\end{document}